\documentclass[a4paper,10pt]{article}

\usepackage{amssymb,amsfonts,amsmath,geometry}
\usepackage[latin1]{inputenc} 
\usepackage[hypertex]{hyperref}

\newtheorem{theorem}{Theorem}
\newtheorem{corollary}{Corollary}
\newtheorem{lemma}{Lemma}
\newtheorem{proposition}{Proposition}
\newtheorem{definition}{Definition}

\newcommand{\nat}{\mathbb{N}}

\newcommand{\re}{\mathbb{R}}

\newcommand{\ber}{\mathcal{B}}

\newcommand{\cqd}{\nopagebreak\hfill\fbox{ }}
\let\a=\alpha

\let\y=\upsilon

\newcommand{\R}{\mathbb{R}}

\newcommand{\N}{\mathbb{N}}

\newcommand{\rar}{\rightarrow}
\newcommand{\LL}{\mathcal{L}}

\newcommand{\TT}{\mathcal{T}_{s,A}}
\newcommand{\LLW}{\mathcal{L}_{A}}
\newcommand{\DLLW}{\mathcal{L}^*_{A}}
\newcommand{\evc}{\psi_A}
\newcommand{\evl}{\lambda_A}
\newcommand{\II}{\int_{M}}
\newcommand{\BB}{\mathcal{B}}

\begin{document}

\title{Entropy and Variational Principle for  one-dimensional Lattice Systems with a general a-priori probability: positive and zero temperature}
\author{A. O. Lopes, J. K. Mengue, J. Mohr and  R. R. Souza} 
\date{\today}

\maketitle

\centerline{Instituto de Matem\'atica, UFRGS - Porto Alegre, Brasil}
\bigskip

We generalize several results of the classical theory of Thermodynamic Formalism  by considering a    compact metric space $M$ as the state space.
We analyze the  shift acting on $M^\mathbb{N}$ and consider a  general a-priori probability for defining the Transfer (Ruelle) operator.
We study potentials $A$ which can depend on the infinite set of coordinates in $M^\mathbb{N}.$
We define entropy and by its very nature it is always a nonpositive number. The concepts of entropy and transfer operator are linked. If M is not a finite set there exist
Gibbs states with arbitrary  negative value of entropy.
Invariant probabilities with support in a fixed point will have entropy equal to minus infinity. In the case $M=S^1$, and the a-priori measure is Lebesgue $dx$, the infinite product of $dx$ on $(S^1)^\mathbb{N}$ will have zero entropy.

We analyze the Pressure problem for a Hölder potential $A$ and its relation with eigenfunctions and eigenprobabilities of the Ruelle operator.
Among other things we analyze the case where temperature goes to zero and we show some selection results.
Our general   setting can be adapted  in order to analyze the Thermodynamic Formalism for the Bernoulli  space with countable infinite symbols. Moreover, the so called $XY$ model also fits under
our setting. In this last case M is the unitary circle $S^1$. We explore the differentiable structure of $(S^1)^\mathbb{N}$ by {\bf considering a certain class of smooth potentials and we show some properties of the corresponding main eigenfunctions}.

\section{Introduction}

Let $(M,d_M)$ be a compact metric space. We consider the metric in $M^{\mathbb{N}}$ given by:
\[d(x,y) = \sum_{n=1}^{\infty}\frac{1}{2^{n}}d_M(x_n,y_n),\]
where $x=(x_1,x_2,...)$ and $y=(y_1,y_2,...)$.
Note that $\ber:=M^{\nat}$ is compact by Tychonoff´s theorem.

We denote by $H_\alpha$ the set of $\alpha$-Hölder functions $A:\ber\to\mathbb{R}$ with the norm
\[\left\|A\right\|_{\alpha} = \|A\| + |A|_{\alpha},\]
where
\[\|A\| = \sup_{x\in\ber}|A(x)| \ \ \text{and} \ \ |A|_{\alpha}=\sup_{x\neq y} \frac{|A(x)-A(y)|}{d(x,y)^{\alpha}}.\]
$\sigma:\ber \to \ber$ denotes the shift map which is defined by
\[\sigma(x_{1},x_{2},x_{3},...)=(x_2,x_3,x_4,...).\]

Let   ${\cal C}$ be the space of continuous functions from $\ber$ to $\R$, and we  will fix an {\it a-priori} probability measure $\nu$ on the Borel sigma algebra over  $M$. We assume that the support of $\nu$ is the set $M$. We stress the crucial point: $\nu$ needs to be a {\it probability measure}, not only a {\it measure}. Note that from our hypothesis if $x_0$ is isolated then $\nu(x_0)>0$. 


For a fixed potential $A \in H_{\alpha}$  we define a Transfer Operator (also called Ruelle operator) $\LL_A: \cal C \to \cal C$ by the rule

$$\LL_A(\varphi)(x) = \int_{M} e^{A(ax)}\varphi(ax) d\nu(a)\,, $$
where  $x\in\ber$ and $ax=(a,x_1,x_2,....)$ denote a pre-image of $x$ with
$a\in M$.
\bigskip

We call One-dimensional Lattice System Theory this general setting. Rigorous mathematical  results on Statistical Mechanics are presented in \cite{Mayer}, \cite{Ge}, \cite{Israel}, \cite{EFS}, \cite{van}, \cite{Si}, \cite{AFMT}, \cite{Be}, \cite{Gro}, \cite{Ellis} and \cite{Gal}.

We point out that a Holder potential $A$ defined on $M^\mathbb{Z}$ is coboundary with a potential in $M^\mathbb{N}$ (same proof as in \cite{PP}).
In this way the Statistical Mechanics of interactions  on $M^\mathbb{Z}$ can be understood via the analysis of the similar problem in $M^\mathbb{N}$.

In \cite{BCLMS} it was investigated the Gibbs measure at positive and zero temperature for a potential $A$ (which depends on infinite coordinates) in the case $M=S^1$, where the a-priori measure is Lebesgue measure. This is the so called  $XY$ model (see \cite{LMST},\cite{Th},\cite{FHo}) which is considered in several applications to real problems in Physics. The spin in each site of the lattice is described by an angle from $[0,2 \pi)$. In the Physics literature, as far as we know, the potential $A$ depends on two coordinates.
A well known example in applications is the potential $A(x)=A(x_0,x_1)= \cos(x_1-x_0-\alpha)+ \gamma\cos(2 x_0)$.

The present paper is a generalization of the setting presented in \cite{BCLMS} for positive and zero temperature.
 We will also consider here a topic which was not addressed there, namely,  the equilibrium (maximizing pressure) measure
for the potential $A$.

There are several possible points of view for understanding Gibbs states  in Statistical Mechanics (see \cite{Sa2}, \cite{Ruelle} for interesting discussions). We prefer the transfer operator method
because we believe that the eigenfunctions and  eigenprobabilities (which can be derived from the theory) allow a more deep understanding of the problem. For example, the information one can get from the main eigenfunction (defined in the whole lattice) is worthwhile, mainly in the limit when temperature goes to zero.

\bigskip

\textbf{Examples:}\newline

 Now we give a brief description of some other examples that fit in our setting.
The last example will be explained in details in section \ref{non}.

\begin{itemize}
\item
 If the alphabet is given by $M=\{1,2,...,d\}$, and the a-priori measure is given by $\displaystyle\nu=\frac{1}{d}\sum_{i=1}^d \delta_i$, then we have the original full shift in a finite set of $d$ symbols and the transfer operator is the classical Ruelle operator associated to a potential $A-\log(d)$ (see for example \cite{PP} and \cite{Kel}). More precisely
\[\LL_A(\varphi)(x) =\int_{M} e^{A(ax)}\varphi(ax) d\nu(a) = \sum_{a\in\{1,2,...,d\}} e^{A(ax)-\log(d)}\varphi(ax).\]
\newline
  If we change the a-priori measure to  $\displaystyle\nu=\sum_{i=1}^d p_i.\delta_i,\,\,\mbox{where }\; p_i>0, \,\,$ and $ \displaystyle\sum_{i=1}^d p_i = 1$, then
\[\LL_A(\varphi)(x) =  \sum_{a\in\{1,2,...,d\}}e^{A(ax)}\varphi(ax)p_a =\sum_{a\in\{1,2,...,d\}}e^{A(ax)+\log(p_a)}\varphi(ax)\]
is the classical Ruelle operator with potential $A+\log(P)$, where $P(x_1,x_2,...)=p_{x_1}$.
\newline
%


\item  If $M_0 = \{z_i, i \in \mathbb{N} \}$ is a countable infinite subset of $S^1$, where each point is isolated, and there is only one accumulating point $z_{\infty} \in S^1 \verb"\" M_0$, then $M=M_0 \cup \{ z_{\infty}\}$
is a compact set. In this case
$M$ can be identified with $\mathbb{N}$, where a special point $z_{\infty}$ plays the role of  infinity (that is, a one-point compactification). We consider here the restricted distance we get from $S^1$ in $M$. If $\displaystyle\sum_{i \in \mathbb{N}} p_i =1 $ with $p_i \geq 0$
and $\displaystyle\nu = \sum_{i \in \mathbb{N}} p_i \delta_{z_i} $ then $\nu$ is supported on the whole $M$, but $z_{\infty}$ is not an atom for $\nu$. The Thermodynamic  Formalism with state space $\mathbb{N}$, or $\mathbb{Z}$,  is considered for example in \cite{Sa},\cite{Sa2},\cite{Da},\cite{JMU},\cite{Mor},\cite{Iom}.  We will analyze in section \ref{non}  some of these results on the present setting.

\end{itemize}


Our main purpose here is to describe a general theory for the Statistical Mechanics of one-dimensional spin lattices.
We point out that most of the papers on the subject assume that the potential $A$ depends just on two (or, a finite number of) coordinates (as for instance is the case of \cite{AFMT},\cite{Be}, \cite{Gro}).
We consider potentials which can depend on the infinite set of coordinates in $M^\mathbb{N}.$

In section \ref{sec:varprinc} we consider the entropy, pressure and Variational Principle and its relations with eigenfunctions and eigenprobabilities of the Ruelle operator. This setting, as far as we know,  was not considered before. In this case the entropy, by its very nature, is always a nonpositive number.  If $M$ is not a  finite set, invariant probabilities with support in a fixed point will have entropy equal to minus infinity. The infinite product of $d\nu$ on $M^\mathbb{N}$ will have zero entropy.
 We point out that, although  at first glance, the fact that the entropy we define here is negative may look  strange, our definition is the natural extension of the concept of Kolomogorov entropy. In the classical case, the entropy is positive because the a-priori measure is not a probability: is the counting measure.

Entropy and Pressure were  considered before in other settings, as for instance in section II in \cite{Israel} or \cite{EFS}. In these works  the authors consider a variational principle on
boxes of finite length, and then they get the equilibrium as the limit probability on the lattice, when the size of the box goes to infinity. The concept of entropy was considered relative to a certain probability on the box (which in some sense plays the role of the a-priori probability). Our formalism is  derived from the Ruelle operator point of view
and is close to the approach described for instance in \cite{PP}, where the probabilities are consider directly on  $M^\mathbb{N}$. As we will see the concepts of entropy and the transfer operator are very much related. When the potential $A$ depends on an infinite number  of coordinates in the lattice we believe our approach is more simple to state and to understand.

Other authors in previous works also considered  Entropy and Transfer Operators on one-dimensional Spin Lattices over metric spaces (see for instance section III \cite{Mayer}, or section A3, or Proposition A4.9 in \cite{Lan}), but we belive our approach is  different.

Among other things we consider in section \ref{zero} the case where temperature goes to zero and show some selection results related with the Ergodic Optimization (see \cite{CGui}, \cite{Jen}, \cite{CLT}, \cite{GL}, \cite{Lep}, \cite{Sou}).
Using the variational principle we obtain a simple proof of the fact that Gibbs states converge to maximizing measures when the temperature goes to zero (a question not discussed in \cite{BCLMS}).


{\bf
An important issue that does not appear in the classical Thermodynamic Formalism (in the sense of \cite{PP} and \cite{Kel}) is the differentiable structure.
We will show in section \ref{dif} that for a certain class of smooth potentials $A$  the associated main eigenfunction is also smooth.}


\section{Ruelle operator}\label{ruellesection}

Let $a^n$ be  an element of $M^n$ having coordinates $a^n = (a_n,a_{n-1},\hdots,a_2,a_1)$, we denote by $a^n x \in \cal{B}$  the concatenation of $a^n \in M^n$ with $x \in {\cal B}$, i.e., $a^n x=(a_n,\hdots,a_1,x_1,x_2,\hdots)$. In the case of $n=1$ we will write $a:=a^1\in M$, and $ax=(a,x_1,x_2,\hdots)$.

The $n$-th iterate of $\LL_A$ has the following expression  $${\cal L}_{A}^n(\varphi)(x)= \int_{M^n}  e^{S_nA(a^n x)} \varphi(a^n x) d\nu^n(a^n),$$ where
$\displaystyle S_nA(a^n x)= \sum_{k=0}^{n-1}A(\sigma^k (a^n x))$ and $\displaystyle d\nu^n(a^n) = \prod_{k=1}^n d\nu(a_{n-k+1})$.

Let us show that $\LL_A$ preserves the set of Hölder functions.

\begin{lemma}  If $\varphi\in H_\alpha$ then $\LL_A(\varphi) \in H_{\alpha}$.
\end{lemma}

\noindent
{\bf Proof:} We have
\[\frac{|\LL_A(\varphi)(x)-\LL_A(\varphi)(y)|}{d(x,y)^{\alpha}} = \frac{|\int_{M}e^{A(ax)}\varphi(ax)d\nu(a) - \int_{M}e^{A(ay)}\varphi(ay)d\nu(a)|}{d(x,y)^{\alpha}}. \]
Now we  use the fact that if $\varphi,A\in H_{\alpha}$, then $e^A \varphi\in H_{\alpha}$, and hence
\[\frac{\int_{M}|e^{A(ax)}\varphi(ax)-e^{A(ay)}\varphi(ay)|d\nu(a)}{d(x,y)^{\alpha}}\leq Hol(e^A \varphi) \nu(M)= Hol(e^A \varphi) .\]

\cqd

\bigskip


\begin{theorem}\label{RPF_eig_nonnorm} Consider a fixed a priori probability $\nu$.
Let us fix $A\in H_{\alpha}$, then there exists a strictly positive Hölder eigenfunction
 $\evc$
for $\LLW: {\cal C} \to {\cal C}$ associated to a strictly positive eigenvalue
$\evl$. This eigenvalue is simple, which means the eigenfunction is unique (modulo multiplication by constant).
\end{theorem}
\textbf{Proof:}
 For each  $0<s<1$, we define the operator $\TT$ on $\mathcal{C}$,  given by
$$\TT(u)(x) = \log\left( \II e^{A(ax)+su(ax)}\, d\nu(a) \right).$$

The introduction of the parameter $s$ in the proof is an adaptation
of an argument presented in \cite{Bousch-Walters}
for the present setting.

{\bf
An easy adaptation of the proof of the proposition 1 in \cite{BCLMS} shows that
$\TT$ is an uniform contraction map. Let $u_s$ be the unique fixed point for $\TT$, then $u_s$ satisfies
\begin{equation}\label{logeigforS}
\log\left(\II e^{A(ax)+su_s(ax)}\, d\nu(a)\right)=u_s(x)\,.
\end{equation}
By the same arguments used in the proof of proposition 2 in \cite{BCLMS}, we can prove that the family $\{ u_s\}_{0<s<1}$ is an equicontinuous family of functions.
It follows from equation \eqref{logeigforS} that
 $$-||A||+ s \min u_s \leq u_s(x)\leq ||A||+ s \max u_s.$$
Hence, $-||A|| \leq (1-s) \min u_s\leq (1-s) \max u_s\leq||A||$, for any $0<s<1$.

 The family $\{ u^*_s=u_s-\max u_s \}_{0<s<1}$ is equicontinuous and uniformly bounded. Let us fix  a subsequence $s_n \to 1$ such that $ [\,(1-s_n) \, \max u_{s_n}\,]\,\to k$, and that, using Arzela-Ascoli theorem,  $\{u_{s_n}^*\}_{n\geq 1}$  has an accumulation point in $\mathcal{C}$, which we will call $u$.

Observe that for any $s$
\begin{eqnarray*}
e^{u^*_s (x)}&=& e^{u_s(x) - \max u_s}=e^{-(1-s) \max u_s + u_s (x) - s \max u_s}\\
&=&e^{-(1-s) \max u_s}\,\int_M e^{  A(ax)+ (s  u_s (ax) - s \max u_s)}\, d\nu(a).
\end{eqnarray*}

Taking limit where $n$ goes to infinity for the sequence $s_n$ we get that $u$ satisfies
$$e^{u(x)} =e^{-k} \,  \int_M e^{  A(ax)+   u (ax)}\,d\nu(a)=e^{-k}\LL_A(e^u)(x).$$

This shows that $\psi_{A}:=e^{u}$ is a positive Holder eigenfunction for $\LL_A$ associated to the eigenvalue $\lambda_A:=e^{k}$.}

The proof of the uniqueness is exactly the same one presented in \cite{BCLMS} (see comments after Theorem 3 in that paper).

\cqd

The eigenfunction is unique up a multiplicative factor. There several ways to normalize it. We assume in this moment that the maximum of the eigenfunction is equal to $1$.
\bigskip

We point out that is possible to generalize the above result for a priori probabilities which depend on the point $x \in M$. This will require some mild assumptions  on this family of probabilities.  We will not address this question here.

\bigskip

We say that a potential $B$ is normalized if $\mathcal{L}_B(1)=1$, which means it satisfies
$$\II e^{B(ax)}d\nu(a) = 1 \,,\,\forall x \in \mathcal{B}\,.$$
In particular, $\forall x \in  \mathcal{B}, a \rar e^{B(ax)}d\nu(a)$ is a probability  measure on $M$,
and $\LL_{B} u(x)$ can be seen as the expectation of the random variable $u$ with respect to this probability measure defined by the point $x$.

Let $A\in H_{\alpha}$,  $\psi_A$ and $\lambda_A$ given by theorem \ref{RPF_eig_nonnorm}, it is easy to see that
\begin{equation}\label{normalizado}\II \frac{e^{A(ax)}\evc(ax)}{\evl \evc(x)}\,\,d\nu(a)=1\,,\,\forall x \in \mathcal{B}\,.\end{equation}
Therefore we define the normalized potential $\bar A$ associated to $A$, as
\begin{equation}\label{pot norm}\bar A := A+\log  \evc-\log \evc \circ \sigma -\log \evl,\end{equation}
where $\sigma:\BB \rar \BB$ is the  shift map.
 As $\psi_A\in H_{\alpha}  $ we have that $\bar A\in H_{\alpha}$. 



We define the Borel sigma-algebra ${\cal F}$ over $\BB$ as  the $\sigma$-algebra
generated by the {\it cylinders}. By this we mean the sigma-algebra generated by sets of the form
$B_1\times B_2\times\,\hdots\,\times B_n \times M^{\mathbb{N}}$, where $n \in \mathbb{N}$, and $B_j, j \in \{1,2,\hdots,n\}$, are open sets in $M$.

We say a probability measure $\mu$ over ${\cal F}$ is invariant, if for any Borel set $B$, we have that $\mu(B)= \mu(\sigma^{-1} (B)).$ 
We denote by ${\cal M}_\sigma$ the set of invariant probability measures. 

We note that $\ber$ is a compact metric space and by the Riesz Representation Theorem, a probability measure on the Borel sigma-algebra is identified with a positive linear functional $L:\mathcal{C} \to \mathbb{R}$ that sends the constant function 1 to the real number 1. We also note that $\mu \in {\cal M}_\sigma$ if and only if, for any $\psi\in \mathcal{C}$ we have
\[ \int_{\ber} \psi \, d\mu = \int_{\ber} \psi\circ \sigma \, d\mu\,.\]

We  define the dual operator $\DLLW$ on the space of Borel measures on $\BB$ as the operator that sends a measure $\mu$ to the measure $\DLLW(\mu)$,  defined by
$$\int_{\BB} \psi\, d\DLLW(\mu) =  \int_{\BB} \LLW(\psi)\, d\mu\,,
$$
for any $\psi \in \mathcal{C}$.


\bigskip

The next theorem is a generalization of propositions 4 and 5 of \cite{BCLMS}.
 Here we consider $\LL_A:H_{\alpha}\to H_{\alpha}$.

\begin{theorem}\label{RPF_normalized}

Let $A$ be a Hölder continuous potential, not necessarily normalized, $\psi_A$ and $\lambda_A$ the eigenfunction and eigenvalue given by  the Theorem \ref{RPF_eig_nonnorm}.
We associate to $A$ the normalized potential
$\bar A = A+\log  \evc-\log \evc \circ \sigma -\log \evl$. Then

(a)  there exists an unique fixed point $\mu_{A}$ for $ {\cal L }_{\bar A}^*$,  which  is a $\sigma$-invariant probability measure;

(b) the measure $$\rho_A= \frac{1}{\psi_A} \,\, \mu_{A}$$ satisfies ${\cal  L}_A^* (\rho_A)= \lambda_A \rho_A$. Therefore, $\,\rho_A$  is an eigen-measure for ${\cal  L}_A^*$;

(c) for any Hölder continuous function $w:{\cal B} \to \mathbb{R}$, we have that, in the uniform convergence topology,
 $${\cal L }_{\bar A}^n\omega \rar \int_{\BB} \omega d\mu_{A}$$
 and  $$\frac{{\cal L}_A^n (w)}{(\lambda_A)^n} \to \, \psi_A \int_{\ber}\,w \,d \rho_A\,,$$
where ${\cal L }_{ A}^n$  denotes the $n$-th iterate of the operator ${\cal L }_{ A}:H_{\alpha}\to H_{\alpha}$.
\end{theorem}

\textbf{Proof:} {\bf
\medskip
(a) We begin by observing that the normalization property implies that the convex and compact set of Borel probability measures on $\BB$ is preserved by the operator ${\cal L}_{\bar A}^*$.
Therefore, using the Tychonoff-Schauder theorem we conclude the existence of a fixed point
$\mu_A$ for the operator ${\cal L}_{\bar A}^*$.
Now we prove that $\mu_A$ is $\sigma$-invariant:  if $\psi \in \mathcal{C}$, we have
$$\int_{\BB} \psi \circ \sigma d\mu_A=\int_{\BB} \psi \circ \sigma d{\cal L}_{\bar A}^*(\mu_A) = \int_{\BB} {\cal L}_{\bar A}( \psi \circ \sigma) d\mu_A =\int_{\BB}  \psi d\mu_A ,
$$
where in the last equality we used the normalization hypothesis for $\bar A$.
The uniqueness of the fixed point will be obtained in the proof of item (c).

(b) ${\cal L}_{\bar A}^*(\mu_A)=\mu_A$ implies that, for any $\psi \in \mathcal{C}$,
\begin{eqnarray*}
  \int_{\ber} \psi d\mu_A &=& \int_{\ber} \psi d {\cal  L}_{\bar A}^*(\mu_A)=\int _{\ber} {\cal  L}_{\bar A}(\psi) d\mu_A \\
   &=& \int_{\ber} \left( \int_M \psi(ax) \frac{e^{A(ax)}  \psi_A(ax)}{\lambda_A \psi_A(x)} d\nu(a)  \right) d\mu_A(x) \,.\\
\end{eqnarray*}
Now, if $\varphi \in \mathcal{C}$, making $\psi=\frac{\varphi}{\psi_A}$ in the last equation, we have
$$\int_{\ber} \frac{\varphi}{\psi_A} d\mu_A = \frac{1}{\lambda_A} \int_{\ber} \left( \int_M \varphi(ax) \frac{e^{A(ax)} }{ \psi_A(x)} d\nu(a)  \right) d\mu_A(x)\,,$$
which is equivalent to
\begin{equation}\label{gibbs_eig}
    \lambda_A \int_{\ber} \varphi d \rho_A = \int_{\ber} {\cal  L}_A(\varphi)  d \rho_A \,,
\end{equation}
i.e., ${\cal  L}_A^* (\rho_A)= \lambda_A \rho_A\,.$

\bigskip

(c) In order to prove item (c) we will need two claims. The first claim can be proved by induction.
\medskip

{\it First Claim:} For any  normalized Holder  potential $B$, if $\|w\|$ denotes the uniform
norm of the Holder function $w:{\cal B}\to \mathbb{R}$,
 we have 
$$|\mathcal{L}_B^n (w) (x) - \mathcal{L}_B^n(w) (y) |\leq \left[C_{e^B}\|w\| \left( \frac{1}{2^{\a}}+...+\frac{1}{2^{n\a}}\right)+\frac{C_w}{2^{n\a}} \right]
d(x,y)^{\a},$$
where $C_{e^B}$ is the Holder constant of $e^B$ and $C_w$ is the Holder constant of $w$.

\medskip

As a consequence of the first claim, the set $\{{\cal  L}_{\bar A}^n \omega\}_{n\geq 0}$ is equicontinuous.
In order to prove that $\{{\cal  L}_{\bar A}^n \omega\}_{n\geq 0}$ is uniformly bounded we  use again the normalization condition, which implies $\|{\cal  L}_{\bar A}^n \omega\|\leq\|w\| \,,\forall n \geq 1$.

Therefore, by the Arzela-Ascoli Theorem there exists  an accumulation point, $\bar \omega$, for  $\{{\cal L}_{\bar A}^n \omega\}_{n\geq 0}$, i.e.,  there exists a subsequence $\{n_k\}_{k\geq 0}$ such that
\begin{equation}\label{wbarra}
\bar \omega (x) = \lim_{k \geq 0} {\cal L}_{\bar A}^{n_k}\omega(x)\,.
\end{equation}

{\it Second Claim:}  $\bar \omega$ is a constant function.

\medskip


To prove this claim, we begin by observing that
\begin{equation}\label{ineqwbarra}
 \sup \bar \omega  \geq \sup {\cal L}_{\bar A} \bar \omega
\end{equation}
 (in fact this inequality holds for any function $w$). Now, \eqref{wbarra} implies
$$ \bar \omega (x) = \lim_{k \geq 0} {\cal L}_{\bar A}^{n_k}\bar \omega(x)\,,$$
(possibly by a different subsequence) and this shows that what we have in \eqref{ineqwbarra} is indeed an equality: in fact, we have
$$  \sup \bar \omega  = \sup {\cal L}^n_{\bar A} \bar \omega \;\forall \; n\geq 0\,.$$
Now, let $x_M^n$ be a maximum point of ${\cal L}^n_{\bar A} \bar \omega$, for any $n \geq 0$.
We have

$$ \bar \omega (x_M^0)= {\cal L}^n_{\bar A} \bar \omega (x_M^n)$$ and this proves the second claim, because the normalization property implies that
$ {\cal L}^n_{\bar A} \bar \omega (x_M^n)$ is a convex combination of $\bar \omega$ in the pre-images of $x_M^n$
(here we also use the fact that  the support of the a-priori probability is the {\it all} space $M$).

Now that $\bar \omega$ is a constant function we can prove  that
$$ \bar \omega = \int_{\BB} \bar \omega d\mu_A = \lim_k \int_{\BB} {\cal L}_{\bar A}^{n_k} \omega d\mu_A
= \lim_k \int_{\BB}  \omega d ({\cal L}_{\bar A}^*)^{n_k}(\mu_A)
= \int_{\BB}  \omega d\mu_A,$$
which shows that $\bar \omega$ does not depend on the subsequence chosen. Therefore, for any $x \in \BB$ we have $${\cal L}_{\bar A}^{n} \omega (x) \rar \bar \omega = \int_{\BB} \omega d\mu_A\,.$$
The last limit shows that the fixed point $\mu_A$ is unique.

To finish the proof of item (c), as
 $A = \bar A - \log  \evc + \log \evc \circ \sigma +\log \evl,$  we have
$$S_{n}A(z) = \sum_{k=0}^{n-1} A\circ \sigma^k(z)=
S_{n}\bar A(z) -\log \psi_A + \log \psi_A \circ \sigma^n + n \log \lambda_A
\,,$$
and therefore
\begin{eqnarray*}
\frac{{\cal L}_A^n(w)(x)}{\lambda_A^n}&=&
\frac{1}{\lambda_A^n}\int_{ M^n} e^{ S_{n}A({a^n}x)} w({ a^n}x) d\nu^n( a^n)=\psi_A(x) \int_{M^n} \frac{e^{ S_{n}\bar A({a^n}x)}}{\psi_A({a^n}x)} w({a^n}x) d\nu^n( a^n)\\
&=&\psi_A(x) {\cal L}_{\bar A}^n\left(\frac{w}{\psi_A}\right)
\to \psi_A(x) \int_{\ber} \frac{w}{\psi_A} d\mu_{A} = \psi_A(x) \int_{\ber} w d\rho_{A} \,.
\end{eqnarray*}}
\cqd

\medskip


We call $\mu_{ A}$ the {\bf Gibbs probability} (or,  Gibbs state) for $A$. We will leave the term {\bf equilibrium probability}  (or,  equilibrium state) for the one which maximizes pressure.
As we will see, this invariant probability measure over ${\cal B}$ describes the statistics in equilibrium for  the interaction described by the potential $A$.
The assumption that the potential is Hölder  implies  that the decay of iteration is fast.


We  normalize $\evc$ by assuming that $\max \evc =1$. There are other possible normalizations.
Therefore, $\rho_A$ is not a probability measure. From the item (c) above we conclude that
\[\frac{{\cal L}_A^n (1)}{(\lambda_A)^n} \to \, \psi_A  \;\rho_A(\ber).\]

\bigskip

\begin{proposition}\label{proponly}

The only Holder continuous eigenfunction $\psi$ of ${\cal L}_A$ which is totally positive is $\psi_A$.

\end{proposition}

{\it Proof:} Suppose $\psi:{\cal B} \to \mathbb{R}$ is a Holder continuous  eigenfunction  of ${\cal L}_A$ associated to some eigenvalue $\lambda$.
It follows from item (c) of theorem \ref{RPF_normalized} that
$$\frac{{\cal L}_A^n (\psi)}{\lambda_A^n} \to \, \psi_A \int_{\ber} \psi d \rho_A\,, \; \mbox{when} \; n \to \infty\;.$$

 Therefore, if $\psi >c>0$, then $\int_{\ber} \psi d \rho_A>0$. Moreover, ${\cal L}_A^n (\psi) = \lambda^ n \psi$. This is only possible if $\lambda=\lambda_A$ and $\psi=\psi_A$.

\cqd

The next result follows from proposition 7 of \cite{BCLMS}.

\begin{proposition}\label{espectro} Suppose $\bar A$ is normalized, then
the eigenvalue $\lambda_{\bar A}=1$ is maximal. Moreover, the remainder of the spectrum of ${\cal L}_{\bar A}: H_\alpha\to H_\alpha$ is contained in a disk centered at zero with radius strictly smaller than one.
\end{proposition}

\vspace{0.3cm}

The proof is the same one as the presented in \cite{PP} (see Theorem 2.2 (ii)).

We denote $\lambda^1_{\bar A}< \lambda_{\bar A}=1$
the spectral radius of ${\cal L}_{\bar A}$ when restricted to the set  $$\{w\in H_\alpha: \int_{\ber} w \,d \mu_{A}=0\}$$
(which is the orthogonal complement of the space of constant functions, i.e., the orthogonal complement of the eigenspace associated to the maximal eigenvalue).
One can also show the exponential decay of correlation for Hölder functions \cite{BCLMS}, which implies mixing and ergodic properties for $\mu_A$.

\begin{proposition}\label{dual e correlacao}
If $v,w \in {\cal L}^2 (\mu_A)$ are such that $w$ is Hölder and $\int_{\ber} w \,  d \mu_{A}=0$, then, there exists $C>0$ such that for all $n$
$$\int_{\ber} (v \circ \sigma^n) \, w \, d \mu_A\leq C \, (\lambda^1_{\bar A})^n.$$
In particular $\mu_A$ is mixing and therefore ergodic.
\end{proposition}

For the proof one can adapt Theorem 2.3 in \cite{PP}.


\section{Entropy and Variational Principle}\label{sec:varprinc}
In this section we will introduce a notion of entropy.  Initially, this will be done only for Gibbs probabilities, and then we will extend this definition to invariant probabilities. After that we prove that the Gibbs probability obtained in the general setting above satisfies a variational principle.
We will also study some general properties of this notion of entropy
and compare it with the classical Kolmogorov entropy when $M=\{1,...,d\}^{\N}$. Finally we will show that this definition is an extension of the notion of entropy for Markov measures (which are the Gibbs measures when the potential depends only on the first two coordinates), as introduced in \cite{LMST}.

Remember that ${\cal M}_\sigma$ is the set of invariant probability measures.

%
%

\begin{definition} Let $\nu$ be a fixed a-priori probability on $M$. We denote by
$ {\mathcal{G}}={\mathcal{G}_\nu}$ the set of Gibbs measures, which means the set of $\mu \in{\cal M}_\sigma $, such that, $\LL^*_{B}(\mu)=\mu$, for some normalized potential $B\in H_{\alpha}$. We define the entropy of $\mu \in \mathcal{G}$ as $$h(\mu)\,=\,h^\nu(\mu)\,=-\int_{\ber} B(x)d\mu(x).$$
\end{definition}

\bigskip

We will see below that $-\int_{\ber} B \, d\mu$ is the infimum of
\[\left\{-\int_{\ber} A\, d\mu + \log(\lambda_A) \,:\, A\in H_{\alpha} \right\}.\]

The above definition is different from the one briefly mentioned in section 3 in \cite{BCLMS}.
\bigskip

\textbf{Remark:} This concept of entropy depends on the choice of the a-priori measure $\nu$, which we choose to be a probability.
In the classical case, when $M=\{1,...,d\}$, the entropy $H(\mu)$ is computed with the a-priori measure $\nu$ given by $\sum_{j=1}^d \delta_j$ (which is not a probability).
A comparison of the value of the above entropy $h(\mu)$, when $M=\{1,...,d\}$, with the classical Kolmogorov entropy $H(\mu)$ (for the full shift) is discussed below, after proposition \ref{entropy-partitions}.
For example, if the a-priori probability is $\nu=\frac{1}{d}\,\sum_{j=1}^d \delta_j$, to get the  entropy $h(\mu)$ you  just have to add $-\log d$ to the classical one $H(\mu)$. Therefore, in this particular case, the above definition results in a number between $-\log(d)$ and $0$. We point out that in the case $M$ has infinite cardinality the above definition $h(\mu)$ makes sense, is well defined, and it is the natural generalization of the previous concept.
\bigskip

\textbf{Remark:} Let $\mu$ be a Gibbs measure and $B$ the normalized potential associated to $\mu$, if we call $J=e^{-B}$,
we have an equivalent definition of entropy given by $$h(\mu)=\int_{\ber} \log (J(x))d\mu(x).$$
We point out that $J=e^{-B}$ does not corresponds to the usual concept of Jacobian of the measure $\mu$.
For example,  consider a finite alphabet $M=\{1,...,d\}$ and $\nu$ the a-priori probability  given by  $\nu(i)=p_i$, where $p_i \geq 0$ and
$\sum_{i=1}^d p_i=1$. In this setting, the Ruelle operator is given by $${\cal L}_B w(x) = \sum_{i=1}^d e^{B(ix)} w(ix) p_i,$$ which can be rewritten as
$${\cal L}_B w(x) = \sum_{i=1}^d e^{B(ix)+\log(p_i)} w(ix)\,.$$ The last formulation fits in the classical thermodynamical formalism setting (see \cite{PP}), for a potential $\tilde B(ix) = B(ix) + \log(p_i)$,  where we know that the Jacobian (defined as $\lim_{n\to \infty} \frac{\mu[x_1x_2...x_n]}{\mu[x_2x_3...x_n]}$, when $[x_1x_2...x_n]$ is the usual cylinder set) is given by $e^{-\tilde B}=e^{-B - \log(p_i)}$. 



\bigskip

\begin{proposition}\label{entropyneg}
If $\mu \in \mathcal{G}$, then we have $h(\mu)\leq 0$.
\end{proposition}

\noindent
{\bf Proof:}
Let $\mu$ be a probability on $\mathcal{G}$ with associated normalized potential $B$. We have
\begin{eqnarray*}
  h(\mu)&=& -\int_{\ber}B(x)d\mu(x)=\int_{\ber}\log e^{-B(x)}d\mu(x) \leq \log \int_{\ber}e^{-B(x)}d\mu(x) \\
   &=& \log \int_{\ber}e^{-B(x)}d\LL^*_{B}(\mu)(x)
 =\log \int_{\ber}\LL_{B}e^{-B(x)} d\mu(x)= 0,
\end{eqnarray*}
where we have used Jensen's inequality and also $\LL_{B}e^{-B(x)}=1$.

\cqd


\bigskip

This negative entropy property will be useful in the next section to give an easy proof that the Gibbs measures of ${\beta A}$ select maximizing measures for $A$, when $\beta\to+\infty$. This result is  obtained in the classical Thermodynamic Formalism Theory because the entropy is bounded. It will be also not difficult  to get this in the present setting because the current notion of entropy is bounded above (by zero).

\medskip

Now we state a lemma that will be used to prove the main result of this section, namely,  the variational principle of Theorem \ref{pvariacional}.
This lemma was shown to be true in  the case $M$ is finite (and the classical Kolmogorov entropy)   in \cite{Lop}.

\begin{lemma}\label{lemavp}
Let us fix a Hölder continuous potential $A$ and a measure $\mu \in \mathcal{G}$ with associated normalized potential $B$.
 We call $\cal C^+$  the space of continuous  positive functions on $\ber$. We have
$$ h(\mu)+\int_{\ber}A(x) d\mu(x)=\inf_{u\in \cal C^+}\bigg\{\int_{\ber} \log \bigg(\frac{\LL_A u(x)}{u(x)}\bigg) d\mu(x)    \bigg\}. $$
\end{lemma}

\noindent
{\bf Proof:} If we take $\tilde u(x)=e^{-A(x)+B(x)}$, then
$$  \log \bigg(\frac{\LL_A \tilde u(x)}{\tilde u(x)}\bigg)=
\log \bigg(\frac{\int_M e^{B(ax)} d\nu(a) }{e^{-A(x)+B(x)}}\bigg) =A(x)-B(x). $$
Integrating, we get
\begin{eqnarray*}
  \int_{\ber} \log \bigg(\frac{\LL_A \tilde u(x)}{\tilde u(x)}\bigg)d\mu(x) &=& \int_{\ber} A(x)d\mu(x)-\int_{\ber}B(x)   d\mu(x) \\
   &=& h(\mu)+  \int_{\ber} A(x)d\mu(x).
\end{eqnarray*}

Now, let us consider a general $\bar u\in \cal C^+$. Using the fact that $e^{-A+B}$ is a positive function, we can write $\bar u(x)=u(x)e^{-A(x)+B(x)}$.
Note that, in this case,
$$ \LL_A\bar u(x)=
\int_M e^{B(ax)} u(ax)  d\nu(a)=\LL_{B}u(x).$$

Hence,
$$  \log \bigg(\frac{\LL_A \bar u(x)}{\bar u(x)}\bigg)= \log( \LL_{B}u(x))-\log u(x)+A(x)-B(x),  $$
and therefore, by integration, we get
\begin{eqnarray*}
 \int_{\ber}   \log \bigg(\frac{\LL_A \bar u(x)}{\bar u(x)}\bigg)d\mu(x)  &=& \int_{\ber} \log( \LL_{B}u(x))d\mu(x)-\int_{\ber}\log u(x)d\mu(x) \\
   &+& \int_{\ber}A(x)d\mu(x)-\int_{\ber}B(x)d\mu(x).
\end{eqnarray*}
Now, all we need to prove is that
$$\int_{\ber} \log( \LL_{B}u(x))d\mu(x)-\int_{\ber}\log u(x)d\mu(x)\geq 0.$$In order to do that,
we use Jensen inequality, and we  get $ \log(\LL_{B} u(x))\geq\LL_{B}\log u(x)$, 
which implies
$$\int_{\ber}\log(\LL_{B} u(x))d\mu(x)  \geq \int_{\ber}\LL_{B}\log u(x)d\mu(x)=\int_{\ber}\log u(x)d\mu(x),$$
where we used $\LL^*_{B}(\mu)=\mu$.



\cqd

\bigskip

Let $\mu \in \mathcal{G}$, with associated normalized potential $B$,  $A\in H_{\alpha}$ and $\psi_A$ and $\lambda_A$, respectively, the positive eigenfunction and the maximal eigenvalue of $\LL_A$, given by theorem \ref{RPF_eig_nonnorm}.
 From Lemma \ref{lemavp} we have
 \begin{eqnarray*}
   h(\mu) + \int_{\ber} A d\mu &=& \inf_{u\in \cal C^+}\bigg\{\int_{\ber} \log  \bigg(\frac{\LL_A u(x)}{u(x)}\bigg) d\mu(x)    \bigg\} \\
    &\leq& \bigg\{\int_{\ber} \log \bigg(\frac{\LL_A \psi_A}{\psi_A}\bigg) d\mu(x)   \bigg\}=\log \lambda_A.
 \end{eqnarray*}
This implies
\[h(\mu)=-\int_{\ber} B d\mu \leq -\int_{\ber} A d\mu+\log \lambda_A,\;\;\;\;  \forall A\in H_{\alpha},\]
with equality if $A=B$ (as $ \lambda_B=1$).
Therefore

\[h(\mu) =  \inf_{A \in H_{\alpha}} \bigg\{-\int_{\ber} A \, d\mu+\log \lambda_A\bigg\},\]
with the minimum  attained at $B$.
Now, based on the last equation, we can extend the definition of entropy for all  invariant measures.

\begin{definition}\label{def:entrinvar}
Let $\mu$ be an invariant measure.
We define the entropy of $\mu$ as
\[h^\nu (\mu)\,=\,h(\mu)\, =\,  \inf_{A \in H_{\alpha}} \bigg\{-\int_{\ber} A \, d\mu+\log \lambda_A\bigg\},\]
where $\lambda_A$ is the maximal eigenvalue of $\LL_A$, given by theorem \ref{RPF_eig_nonnorm}.
\end{definition}

This value is non positive and can be $-\infty$ as we will se later.

\medskip

\begin{definition}\label{def:press} Given a Hölder potential $A$ we call the pressure of $A$ the value

$$  P(A)\,=\,\sup_{\mu\in\cal{M}_{\sigma}} \bigg\{h(\mu)+\int_{\ber}A(x) d\mu(x)\bigg\}.  $$

A probability which attains such maximum value is called  equilibrium state for $A$.

\end{definition}

In the literature sometimes this value is called the Free Energy (see \cite{Ellis}  for instance).

Now we will show  the variational principle of pressure which characterizes the equilibrium state:

\begin{theorem}[Variational Principle] \label{pvariacional}
Let $A\in H_{\alpha}$ be a Hölder continuous potential and $\lambda_A$ be  the maximal eigenvalue of $\LL_A$, then
$$  \log \lambda_{A}\,=\,P(A)\,=\,\sup_{\mu\in\cal{M}_{\sigma}} \bigg\{h(\mu)+\int_{\ber}A(x) d\mu(x)\bigg\}.  $$
Moreover the supremum is attained on the Gibbs measure, i.e. the measure $\mu_A$ that satisfies $\LL^*_{\bar{A}}(\mu_A)=\mu_A$.
\end{theorem}
\medskip

\smallskip

Therefore, the Gibbs state and the equilibrium state for $A$ are given by the same measure $\mu_A$, which is the unique fixed point for the
dual Ruelle operator associated to the normalized potential $\bar{A}$.
\bigskip

\noindent
 {\bf Proof:} Consider a fixed  $A\in H_{\alpha}$, by the definition of entropy, we have
 \begin{eqnarray*}
    & & \sup_{\mu\in\cal{M}_{\sigma}} \bigg\{h(\mu)+\int_{\ber}A(x) d\mu(x)\bigg\} \\
    & & =\sup_{\mu\in\cal{M}_{\sigma}} \bigg\{\inf_{B \in H_{\alpha}} \bigg\{-\int_{\ber} B \, d\mu+\log \lambda_B\bigg\}+\int_{\ber}A(x) d\mu(x)\bigg\} \\
    & & \leq \sup_{\mu\in\cal{M}_{\sigma}} \bigg\{ -\int_{\ber} A \, d\mu+\log \lambda_A+\int_{\ber}A(x) d\mu(x)\bigg\} = \log \lambda_A.
 \end{eqnarray*}
 Hence,  $$\log \lambda_A\geq \sup_{\mu\in\cal{M}_{\sigma}} \bigg\{h(\mu)+\int_{\ber}A(x) d\mu(x)\bigg\}.        $$
 On the order hand, as $A\in H_{\alpha}$, from theorem \ref{RPF_eig_nonnorm}  we know that there exists $\lambda_A$ and $\varphi_A$, such that, $\LL_A(\varphi_A)= \lambda_A\varphi_A$. Now, if we define $\bar A= A+\log\varphi_A-\log\varphi_A\circ \sigma-\log\lambda_A$, then,  by Theorem \ref{RPF_normalized},  there exists a measure $\mu_A$ such that $\LL^*_{\bar A}(\mu_A)=\mu_A$.
 This implies
  $\mu_A\in \mathcal {G}$, and $$h(\mu_A)=-\int_{\ber}\bar A(x)d\mu_{ A}(x)= -\int_{\ber} A \, d\mu_{A}+\log \lambda_A.$$
 Therefore,
 $$ \log \lambda_A =h(\mu_A)+\int_{\ber} A \, d\mu_{A}\leq \sup_{\mu\in\cal{M}_{\sigma}} \bigg\{h(\mu)+\int_{\ber}A(x) d\mu(x)\bigg\}.  $$

\cqd

\bigskip

In \cite{LMST} a  variational principle of pressure was considered. Other variational principles of pressure were described in \cite{Israel} \cite{Mayer}. Our approach and also the kind of probabilities we consider are different of the ones in  this last reference.

\begin{theorem}[Pressure as Minimax] \label{minimax} Given a Hölder potential $A$

$$ P(A)=\sup_{\mu\in\cal{M}_{\sigma}}\,\Big[\,\inf_{u\in \cal C^+}\bigg\{\int_{\ber} \log \bigg(\frac{\LL_A u(x)}{u(x)}\bigg) d\mu(x)    \bigg\}\,\Big]. $$

\end{theorem}

{\bf Proof:} This follows at once from Lemma \ref{lemavp} (see also \cite{LO}).

\cqd

\bigskip


 It is known that periodic orbits can be used to get information about the pressure in the classical thermodynamic formalism setting,  and also to approximate the equilibrium measure (see \cite{PP} chapter 5 and \cite{LM3}).
Therefore, the next corollary can be useful:

\begin{corollary} For each  $a^{n}\in M^{n}$, $a^{n}=(a_n,...,a_1)$  let $a^{\infty}\in \ber$ be the periodic orbit of period $n$ obtained by the successive concatenation of $a^{n}$, i.e., $a^{\infty}=(a_n,...,a_1,a_n,...,a_1,...)$. Then,
\[
P(A)=\lim_{n\to\infty}\frac{1}{n}\log\left(\int_{M^{n}}e^{S_nA (a^{\infty})}\, d\nu^{n}(a^{n}) \right).\]
\end{corollary}

\noindent
{\bf Proof:}
Using Theorem \ref{RPF_normalized} (c) and then the Variational Principle we conclude that for any fixed $x\in \ber$
\[\lim_{n\to\infty}\frac{1}{n}\log\left(\int_{M^{n}}e^{S_nA (a^{n}x)}\, d\nu^{n}(a^{n}) \right)=
\lim_{n\to\infty}\frac{1}{n}\log\left(\LL_A^n(1)(x) \right)
= \log(\lambda_A)=
P(A).\]
Using the fact that $A$ is Hölder continuous, there exists a constant $C>0$ such that
\[|S_nA (a^{n}x) - S_nA(a^{\infty})|\leq C\left(\frac{1}{2^{\alpha}}+\frac{1}{2^{2\alpha}}+...+\frac{1}{2^{n\alpha}}\right)\,d(x,a^{\infty})^{\alpha}.\]
Therefore, using   that $\ber$ is a compact set with finite diameter, there exists a constant $C>0$, such that,
\[|S_nA (a^{n}x) - S_nA(a^{\infty})|\leq C, \]
for any $n\in\mathbb{N},x\in\ber$ and $a\in M$. Then,
\[\int_{M^{n}}e^{S_nA (a^{n}x)-C}\, d\nu^{n}(a^{n})\leq \int_{M^{n}}e^{S_nA (a^{\infty})}\, d\nu^{n}(a^{n})\leq \int_{M^{n}}e^{S_nA (a^{n}x)+C}\, d\nu^{n}(a^{n}),\]
and
\[\lim_{n\to\infty}\frac{1}{n}\log\left(\int_{M^{n}}e^{S_nA (a^{\infty})}\, d\nu^{n}(a^{n}) \right)=\lim_{n\to\infty}\frac{1}{n}\log\left(\int_{M^{n}}e^{S_nA (a^{n}x)}\, d\nu^{n}(a^{n}) \right)= P(A).\]
\cqd

\bigskip

Now we present a few  properties of entropy:

\begin{proposition}\label{def:entrinvar1}The entropy has the following properties:
\newline a) $h(\mu) \leq 0$ for any invariant measure $\mu$.\newline
b) $\nu^{\infty}=\nu\times\nu\times\nu\times...$ has zero entropy.\newline
c) The entropy is upper semi-continuous. \newline
d) The entropy is a concave function in the space of invariant probabilities. \newline
\end{proposition}

\noindent
{\bf Proof:}
\begin{itemize}
  \item[a)]  We point out that $A=0$ is a normalized function, hence for any invariant measure $\mu$  we have $h(\mu) \leq 0$.\newline
  \item[b)] We are going to show that $\nu^{\infty}$ is the equilibrium measure for $A=0$. Indeed,
      \begin{eqnarray*}
        \int_{\ber} g \, d\LL^{*}_{A}(\nu^{\infty}) &=& \int_{\ber} \LL_A g\, d\nu^{\infty} = \int_{\ber} \int_{M} e^{0} g(ax)\,d\nu(a)d\nu^{\infty}(x) \\
         &\stackrel{z=ax}{=}& \int_{\ber}g(z)\, d\nu^{\infty}(z).
      \end{eqnarray*}
   which shows that $\LL^{*}_{A}(\nu^{\infty})=\nu^{\infty}$.
  \item[c)] Fix an $\varepsilon>0$ and suppose $\mu_n$ converges to $\mu$.
By definition of $h(\mu)$,  we can choose  $A\in H_{\alpha}$  such that
$$-\int_{\ber} A d\mu + \log \lambda_A \leq h(\mu) + \varepsilon.$$
If $n$ is large enough, we have
$| \int_{\ber} A d\mu_n -\int_{\ber} A d\mu| < \varepsilon$.
Then, 
$$h(\mu_n)  \leq -\int_{\ber} A d\mu_n +\log \lambda_A \leq -\int_{\ber} A d\mu + \log \lambda_A +\varepsilon   \leq   h(\mu)+2\varepsilon\,,$$
therefore
$\displaystyle{\limsup_{n\to+\infty} h(\mu_n) \leq h(\mu)+2\varepsilon}.$
  \item[d)] Let $\mu_1$ and $\mu_2$ be $\sigma$-invariant probabilities, $\varepsilon \in (0,1)$ and $\mu=\varepsilon\mu_1 + (1-\varepsilon)\mu_2$. Then
   \begin{eqnarray*}
      & & h(\varepsilon\mu_1 + (1-\varepsilon)\mu_2) = h(\mu) = \inf_A\left(-\int_{\ber} A \, d\mu + \log\lambda_A\right) \\
      & & =\inf_A\left(-\varepsilon \int_{\ber} A \, d\mu_1+(1-\varepsilon)\int_{\ber} A \, d\mu_2 + \log\lambda_A\right) \\
      & & \geq \inf_A\left(-\varepsilon \int_{\ber} A d\mu_1+\varepsilon.\log\lambda_A\right)+\inf_A\left(-(1-\varepsilon) \int_{\ber} A d\mu_2+(1-\varepsilon).\log\lambda_A\right) \\
      & & =\varepsilon h(\mu_1)+(1-\varepsilon)h(\mu_2).
   \end{eqnarray*}
 \end{itemize}
\cqd


\bigskip

\textbf{Remark:} The entropy of a probability measure supported on periodic orbit can be $-\infty$.
%
%
Indeed, suppose $M=[0,1]$, and $A_c:M^{\nat}\rar \re$  given by $A_c(x)= \log\left(\frac{c}{1-e^{-c}}e^{-cx_1}\right) $. Suppose the a-priori $\nu$ measure is the Lebesgue measure. We have that for each $c>0$, the function $A_c$ is a $C^1$ normalized potential (therefore belongs to $H_{\alpha}$), which depends only on the first coordinate of $x$. Note that $\mathcal{L}_{A_c}(1)=1$. Let $\mu$ be the Dirac Measure on $0^{\infty}$. We have $h(\mu)\leq -\int_{\ber} A_c d\mu = -A_c(0^{\infty})= - \log\left(\frac{c}{1-e^{-c}}\right) \to -\infty $ when $c \to \infty$.
This shows that $h(\mu)=-\infty$. An easy adaptation of the arguments can be done to prove that, in this setting, 
invariant measures supported on periodic orbits have entropy $-\infty$.
%
%


\bigskip

Note the subtle point  that the entropy depends on the a-priori probability and moreover all subsequent concepts we introduced, like for example the Ruelle operator,
$$\LL_A(\varphi)(x) = \int_{M} e^{A(ax)}\varphi(ax) d\nu(a)\, $$ assume conditions  on the  pre-images of $\sigma$. Therefore, given an iterate $\sigma^n$, if one wants to consider the entropy of a $\sigma^n$-invariant probability, then we need to specify a certain  a-priori probability. We will address this question now.

\bigskip

{\bf Entropy of iterates:}
Suppose $M^n$ is the compact set given by $$\{(x_1,x_2,...,x_n)\, \,| \,x_i \in M \;, \forall \,\,1 \leq i \leq n \}\,,$$ with the sum or the maximum norm. 
Let $\sigma^n$ be the shift map defined on the Bernoulli space given by $\ber^n \equiv (M^n)^\nat$. We know that $\sigma^n$ is the n-th iterate of $\sigma$ in the original Bernoulli space, but we prefer to see $\sigma^n$ as a new map defined on a new Bernoulli space. If we do that, all the theory developed above applies to $\sigma^n$, we have a Ruelle operator with an a-priori measure given by $\nu^n$,  and therefore  the entropy of a Gibbs measure to the new map $\sigma^n$ is well defined. Note that the new Bernoulli set $\ber^n$ can be identified with the original
$\ber$ and an invariant measure for $\sigma$ is also an invariant measure for $\sigma^n$.

\begin{proposition}
 If $\mathcal{G}^n$ denotes the set of Gibbs measures on $\ber^n$,  then $\mu\in \mathcal{G}$ implies $\mu\in \mathcal{G}^n$ and
           \begin{equation}\label{entropia-iterado}
       h_{\mu}^{\nu^n}(\sigma^n)\,= \,    h_{\mu}(\sigma^n)\,=\,n h_{\mu}(\sigma)\,=\, n h_{\mu}^\nu(\sigma)\,.
           \end{equation}

 \end{proposition}

\noindent
{\bf Proof:}
Note that if $\mu$ is a $\sigma$-invariant measure, then $\mu$ is invariant for  $\sigma^n$. Also, if $B \in H_{\alpha}$ is a normalized potential for $\sigma$, then the Birkhoff sum $B^n \equiv \displaystyle\sum_{j=0}^{n-1} B \circ \sigma^j$ is a normalized potential for the map $\sigma^n$.
Let us first prove that if $\mu$ is the Gibbs measure for the  Ruelle operator associated to $B$,
  then, $\mu$ ( which is indeed a measure on $\mathcal{B}^n$), is also Gibbs for $B^n$. In order to do that, note that

  \begin{eqnarray*}
     & & {\cal L}_{B^n}^m(\varphi)(x)= \int_{(M^n)^m}  e^{S_m B^n(a^{nm} x)} \varphi(a^{nm} x) d(\nu^n)^m(a^{nm}) \\
     & &= \int_{(M^n)^m}  e^{S_{mn} B(a^{nm} x)} \varphi(a^{nm} x) d\nu^{nm}(a^{nm}) = {\cal L}_{B}^{mn}(\varphi)(x) \to \int_{\ber} \varphi d\mu.
  \end{eqnarray*}


Now the $\sigma^n$-entropy of $\mu$, given by the integral of the $B^n$, equals $n$ times the $\sigma$-entropy of $\mu$, because, using the fact that $\mu$ is $\sigma$-invariant, we have
 $$ h_{\mu}(\sigma^n) = -\int_{\ber^n} B^n d\mu = -\int_{\ber^n}   \sum_{j=0}^{n-1} B \circ \sigma^j d\mu
 $$
 $$ =    - \sum_{j=0}^{n-1} \int_{\ber} B \circ \sigma^j d\mu = - n \int_{\ber} B d\mu = n h_{\mu}(\sigma)
 \,.$$

\cqd

\bigskip

{\bf Relations with Kolmogorov Entropy:}

\bigskip

Let us consider the construction of the entropy by partitions method, in the case $M$ is finite. We begin by remembering that,  by  the Kolmogorov-Sinai Theorem, the classical entropy of $\mu$, which we will denote by $H(\mu)$, is given by
\begin{equation}\label{KS}
H(\mu) = \lim_{n\to\infty} -\frac{1}{n}\sum_{i_1,...,i_n}\mu([i_1...i_n])\log\left(\mu([i_1...i_n])\right).
\end{equation}

\begin{proposition}\label{entropy-partitions} Let $M=\{1,...,d\}$ and $\displaystyle\nu = \sum_{i=1}^{d}p_i\delta_i$ be the a-priori probability on $M$.  For any Gibbs measure $\mu$:

\begin{enumerate}
\item[(a)]
\[H(\mu) =  h^{\nu}(\mu) - \sum_{i=1}^{d}\log(p_i).\mu([i]),\]

\item[(b)]
\[h^{\nu}(\mu) =- \lim_{n\to\infty} \frac{1}{n}\sum_{i_1,...,i_n}\mu([i_1...i_n])\log\left(\frac{\mu([i_1...i_n])}{p_{i_1}...p_{i_n}}\right)\]
where
\[ [i_1...i_n]=\{x\in M^{\mathbb{N}} : x_1=i_1,...,x_n=i_n\}.\]

\end{enumerate}

\end{proposition}

{\bf Proof:}

If $\mu$ is a Gibbs measure, there exists a  normalized potential $A$ associated to $\mu$, which implies
\[\int_M e^{A(ax)} d\nu(a) = \sum_{i=1}^{d}e^{A(ix)}p_i = 1, \ \ \ \ \forall \, x\in M^{\N},\]
which is equivalent to
\[\sum_{i=1}^{d}e^{A(ix)+\log(p_i)} = 1, \ \ \ \ \forall \, x\in M^{\N}.\]
Moreover, $\LL^*_A(\mu)=\mu$ implies that
$\mu$ is a fixed point  for the Classical Ruelle Operator with the normalized potential $A+\log(P),$ where $ \, P(y_1,y_2,...)=p_{y_1}$.
Therefore
\[H(\mu) = -\int_{\ber} A+\log(P) d\mu = h^{\nu}(\mu) - \int_{\ber} \log(P)d\mu = h^{\nu}(\mu) - \sum_{i=1}^{d}\log(p_i).\mu([i]).\]
which ends the proof of item (a)

In order to prove item (b), we note that, from the last equation, and using that $\mu$ is an $\sigma$-invariant measure, we have, for any $n\geq 1$
\begin{eqnarray*}
  H(\mu) &=& h^{\nu}(\mu)
   - \frac{1}{n}\int_{\ber}\log(P)+...+\log(P\circ\sigma^{n-1})\, d\mu \\
   &=& h^{\nu}(\mu) - \frac{1}{n}\int_{\ber} \log(p_{x_1}...p_{x_n})\, d\mu(x) \\
   &=& h^{\nu}(\mu) - \frac{1}{n}\sum_{i_1,...,i_n}\mu([i_1...i_n])\log(p_{i_1}...p_{i_n}).
\end{eqnarray*}
Then,
\begin{eqnarray*}
  h^{\nu}(\mu) &=& H(\mu) + \frac{1}{n}\sum_{i_1,...,i_n}\mu([i_1...i_n])\log(p_{i_1}...p_{i_n}) \\
   &=& H(\mu) + \lim_{n\to\infty}\frac{1}{n}\sum_{i_1,...,i_n}\mu([i_1...i_n])\log(p_{i_1}...p_{i_n}) \\
   &=& - \lim_{n\to\infty} \frac{1}{n}\sum_{i_1,...,i_n}\mu([i_1...i_n])\log\left(\frac{\mu([i_1...i_n])}{p_{i_1}...p_{i_n}}\right). \end{eqnarray*}

where in the last equation we used \eqref{KS}.


\cqd


In particular, it follows from item (a) above that, when $p_i=\frac{1}{d}$, for all $i$, we have
$$h^{\nu}(\mu) = H(\mu)-\log(d)\,.$$
The above proposition can be interpreted in the following way:  in the classical definition of Kolmogorov entropy it is considered the a-priori measure $\nu = \sum_{i=1}^{\infty}\delta_i$ on $M$, which  is not a probability.









\bigskip

{\bf Markov Chains with values on $S^1$:}

\bigskip


Now we recall the concept of Markov measures and show that the entropy defined above is an extension of the concept of entropy for Markov measures, as introduced in \cite{LMST}.



Let $K:M^2\to \re$, $\theta:M\to \re$, satisfying
\begin{equation}\label {K theta}\int_{M} K(x_1,x_2)d\nu(x_2)=1, \;\forall x_1\;\;\mbox{ and }\;\int_M \theta(x_1)K(x_1,x_2)d\nu(x_1)=\theta(x_2)\;, \forall x_2\,.\end{equation}
We call $K$ a transition kernel and $\theta$ the stationary measure for $K$.
As in \cite{LMST},  we define the absolutely continuous Markov measure associated to $K$ and $\theta$, as

\begin{equation}\label{def:markov}
    \mu(A_1...A_n\times M^\mathbb{N}):=\int_{A_1...A_n}\,\theta(x_1)\,
K(x_{1},x_2)...K(x_{n-1},x_n)\,d\nu(x_n)...d\nu(x_1),
\end{equation}

for any cylinder $ A_1...A_n
\times M^\mathbb{N}$.

The next proposition show us the importance of a.c. Markov measures:

\begin{proposition}\label{MarkovequalGibbs} We will show that
\begin{itemize}
\item[a)]  Given a Hölder continuous potential  $A(x_1,x_2)$ (not necessarily normalized) depending on two coordinates, there exists a Markov measure that is Gibbs for $A$.
\item[b)] The converse is also true: given an absolutely continuous Markov measure defined by $K$ and $\theta$, there exists a certain Hölder continuous normalized potential  $A(x_1,x_2)$, such that the Markov measure defined by  $\theta$ and $K$ is the Gibbs measure for $A$.
\end{itemize}
\end{proposition}

Therefore, any a.c. Markov measure  is Gibbs for a potential depending on two variables, and conversely, any potential depending on two variables has a Gibbs measure which is an a.c. Markov Measure.

In other words, if we restrict our analysis to potentials that depend just on the first two coordinates, we have that the set of a.c. Markov Measures coincides with the set of Gibbs measures. 

\noindent
{\bf Proof:}

(a) Given a potential $A(x_1,x_2)$, non-normalized, as in \cite{LMST} define
$\theta_{A}:M\to \mathbb{R}$ by \begin{equation}\label{theta}\theta_{A}(x_1):=
\frac{\psi_{A}(x_1)\,\, \bar\psi_{A}(x_1)}{\pi_{A}},\end{equation}
and a transition $ K_{A}:M^2\to \mathbb{R}$ by
\begin{equation}\label{K}K_{A}(x_1,x_2):=\frac{e^{ A(x_1,x_2)}\,\,\bar \psi_{A}(x_2)}{\bar \psi_{A}(x_1)\,\lambda_{A}}\,, \;\end{equation}
where
$\psi_A$ and $\bar \psi_A$ are the eigenfunctions associated to the maximal eigenvalue $\lambda_A$ of the operators
\begin{equation}\label{op2coord}    L_A \psi (x_2) =\int_M e^{ A(x_1,x_2)} \, \psi(x_1)d\nu(x_1) \ \ \  \mbox{ and }\ \ \ \bar{ L}_A \psi (x_1) =\int_M e^{ A(x_1,x_2)} \, \psi(x_2)d\nu(x_2)\end{equation}
and   $\pi_{A}=\int_M \psi_{A}(x_1) \bar\psi_{A}(x_1)d\nu(x_1)$.

Then, by the same arguments used to prove theorem 16 of \cite{BCLMS}, we obtain that the Markov measure $\mu_A$ defined by \eqref{def:markov} (considering $K_A$ and $\theta_A$) is Gibbs for $A$, i.e. a fixed point for the dual Ruelle operator $ {\cal L}_{\bar A}^*$, where $\bar A= A+\log  \psi_{A }(x_1)-\log \psi_{A}(x_2)-\log\lambda_A$.

\bigskip

(b) Let $K$ and $\theta$ satisfying \eqref{K theta}, and define $A=\log K$,
we have $\bar{ L}_A(1)=1$ which implies $\lambda_A=1$ and $\bar \psi_A=1$. Let $\psi_A$ be maximal eigenfunction for $L_A$.

 Using \eqref{K}, we get $K_A(x_1,x_2)=e^{A(x_1,x_2)}=K(x_1,x_2)$.
Define $\theta_A=\frac{\psi_A}{\pi_A}$. We have that $\theta_A$ is an invariant density for $K$, therefore $\theta_A=\theta$.
Then, also by theorem 16 page of \cite{BCLMS}, we have that the Markov measure defined by $K$ and $\theta$ is Gibbs for $A$.
\cqd

Next proposition shows that the concept of entropy introduced in \ref{def:entrinvar} is a generalization of the concept of entropy defined in \cite{LMST}, which could only be applied to a.c. Markov measures:

\begin{proposition} \label{twoentr}
Let $\mu$ be the Markov measure defined by a transition kernel  $K$ and a stationary measure $\theta$, given in \eqref{def:markov}.
The definition of entropy given in \cite{LMST}:
$$S(\theta K )= - \int_{M^2} \theta(x_1) K(x_1,x_2) \log(K(x_1,x_2)) d\nu(x_1) d\nu(x_2)\leq 0$$
coincides with the present definition \ref{def:entrinvar}.
\end{proposition}

\noindent
{\bf Proof:} As in  the proof of proposition \ref{MarkovequalGibbs}, note that
the normalized potential associated to  $A(x,y)=\log K(x,y)$ is $$\bar A(x_1,x_2) =\log K(x_1,x_2) +\log  \psi_{A }(x_1)-\log \psi_{A }(x_2), $$
where $\psi_{A }$ is the maximal eigenfunction of the operator $L_A$. Note also that $\bar A$ depends only on the first two coordinates.

Let  $\mu$ be the Gibbs measure associated to $ \bar A$, hence by definition \ref{def:entrinvar} we have
\begin{eqnarray*}
  h(\mu) &=& -\int_\ber \bar A(x_1,x_2) d\mu(x) = -\int_\ber \log K(x_1,x_2) d\mu(x) \\
   &=& \int_{M^2} \log K(x_1,x_2) \theta(x_1) K(x_1,x_2) d\nu(x_1) d\nu(x_2) =S(\theta K).
\end{eqnarray*}
\cqd


\bigskip

\section{Zero temperature}\label{zero}

Consider a fixed Hölder potential $A$ and a real variable $\beta>0$.
We denote, respectively, by $\psi_{\beta A}$ and $\mu_{\beta A}$, the eigenfunction for the Ruelle operator associated to $\beta A$ and the equilibrium measure (Gibbs) for $\beta A$. We would like to investigate general properties of the limits of $\mu_{\beta_n A} $ and of $\frac{1}{\beta_n}\log\psi_{\beta_n A}$ when $\beta_n\to \infty. $ Some results of this section are generalizations of the ones in \cite{BCLMS}. It is well known that the parameter $\beta$ represents the inverse of the temperature.

It is fair to call ``Gibbs state at zero temperature  for the potential $A$" any of the weak limits of convergent subsequences $\mu_{\beta_n A} $. Even when the potential $A$ is Hölder,  Gibbs state at zero temperature do not have to be unique. In the case there exist the weak limit $\mu_{\beta A}\rightharpoonup \mu $, $\beta \to \infty$,
we say that there exists selection of Gibbs state for $A$ at temperature zero.

\bigskip

{\bf Remark :}
Given $\beta$ and $A$, the Hölder constant of $u_{\beta A}=\log(\psi_{\beta A})$, depends on the Hölder constant for $\beta\, A$, and is given by $\beta \frac{2^{\alpha}}{2^{\alpha}-1} Hol_{ A}$
(see \cite{BCLMS}).
As we normalize $\psi_{\beta A}$ assuming that $\max \psi_{\beta A} = 1$, the family of functions $\frac{1}{\beta}\log(\psi_{\beta A}),\, \beta>0$, is  uniformly bounded.
Note that when we normalize $\psi_{\beta A}$ the Hölder  constant of $\log (\psi_{\beta A})$ remains unchanged, which assures the family $\frac{1}{\beta} \log (\psi_{\beta A})\,,\beta>0$, is equicontinuous.

Therefore, there exists a subsequence $\beta_n\to \infty$, and $V$ Holder, such that, on the uniform convergence topology
$$V:=\lim_{n\to\infty}
\frac{1}{\beta_n}\log(\psi_{\beta_n A}).$$




Remember that we denote by $ \mathcal M_\sigma$ the set of $\sigma$ invariant Borel probability measures over ${\cal  B}$. As $ \mathcal M_\sigma$ is compact, given $A$, there always exists a subsequence $\beta_n$, such that $\mu_{\beta_n A}$ converges to an invariant probability measure.

The limits of $\mu_{\beta A}$ are related (see below) with the following problem: given $A:{\cal B} \to \mathbb{R}$ Holder, we want to find probabilities that  maximize, over $ \mathcal M_\sigma$, the value

$$ \int_{\ber} A(x) \,d \mu(\mathbf {x}). $$
We define  
$$m(A)=\max_{\mu\in\mathcal M_\sigma} \left\{ \int_{\ber} A d\mu \right\}\,.$$

Any of the probability measures which attains the maximal value will be called a maximizing probability measure, which will be sometimes denoted generically by $\mu_\infty$. As $\mathcal M_\sigma$ is compact, there exist always at least one maximizing probability measure. It is also true that there
exists ergodic maximizing probability measures. Indeed, the set of maximizing probability measures is convex, compact and the extreme probability measures of this convex set are ergodic (can not be expressed as convex combination of others \cite{Kel}).  Any maximizing probability measure is a convex combination of ergodic ones \cite{PY}.  Results obtained in this setting belong to what is called Ergodic Optimization Theory \cite{Jen}.

The possible limits of $\frac{1}{\beta_n}\log\psi_{\beta_n A}$ are related (see below) with the following concept:

\begin{definition}\label{sub}
A continuous function $u: {\cal B} \to \mathbb{R} $ is called a
{\em calibrated subaction} for $A:{\cal B}\to \mathbb{R}$, if, for any $y\in {\cal B}$, we have

\begin{equation}\label{c} u(y)=\max_{\sigma(x)=y} [A(x)+ u(x)-m(A)].\end{equation}

\end{definition}

This can also be expressed as
$$m(A)= \max_{a \in M} \{A(ay)+  u (ay) - u(y) \} .$$
Note that for any $x\in {\cal B}$ we have
 $$ u(\sigma(x)) -  u(x) - A(x) + m(A) \geq0.$$

The above equation for $u$ can be seen as a kind of discrete version of the concept of sub-solution of the Hamilton-Jacobi equation \cite{CI} \cite{BC} \cite{Fathi}. It can be also seen as a kind of dynamic additive eigenvalue problem \cite{CD} \cite{CG} \cite{GL3}.

We note that $m(A)$ can be characterized by
\[m(A) = \inf \{\gamma : \exists \, u\in {C}, \, \gamma + u\circ\sigma -u -A \geq 0\},\]
where $\cal{C}$ denotes the set of continuous real-valued functions.
In some sense this corresponds to the dual problem in transport theory \cite{LM2}. Any invariant measure help us to estimate $m(A)$ from  below. The continuous functions on the dual problem help us to estimate $m(A)$ from above.

If $u$ is a calibrated subaction, then $u+c$, where $c$ is a constant, is also a calibrated subaction. An interesting question is when such calibrated subaction $u$ is unique up to an additive constant (see \cite{LMST} and \cite{GL}).

Remember that if $\mu$ is $\sigma$-invariant, then for any continuous function $u:{\cal B}\to \mathbb{R}$ we have
$$ \int_{\ber} \,[ u(\sigma(x)) -  u(x)]\, d \mu=0.$$
Therefore if $\mu_{\infty}$ is a maximizing probability measure for $A$ and $u$ is a calibrated subaction for $A$, then
(see for instance \cite{CLT} \cite{Jen} \cite{Sou} for a similar result) for any $x$ in the support of $\mu_\infty$, we have
\begin{equation}\label{funcaoR} u(\sigma(x)) -  u(x) - A(x) + m(A)=0. \end{equation}
In this way if we know the value $m(A)$, then a calibrated subaction $u$ for $A$ can help us to identify the support of maximizing probabilities. The above equation can be eventually true outside the union of the supports of the  maximizing probabilities (see an interesting example due to R. Leplaideur  \cite{CLO}).

We show below that if there exists a subsequence $\beta_n\to \infty$,  such that on the uniform convergence
$$V:=\lim_{n\to\infty}
\frac{1}{\beta_n}\log(\psi_{\beta_n A}),$$
then such $V$ is a calibrated subaction for $A$.
When there exists a $V$ which is the limit
$$V:=\lim_{\beta\to\infty}
\frac{1}{\beta}\log(\psi_{\beta A}),$$
(not just via a subsequence) we say we have selection of subaction at temperature zero. Positive results in this direction are presented in \cite{BLM}, \cite{LM} and \cite{LMST}.

There exists here a subtle point. Sub-action is a concept in Ergodic Optimization and does not depend on the existence of an a-priori probability $\nu$ in $M$. On the other hand, the eigenfunction $\psi_{\beta_n A}$ is associated to a Ruelle Operator, which   depends on the a-priori measure. In any  case, for any given  a-priori probability  $\nu$, if the associated  family of eigenfunctions  $\psi_{\beta_n A}$ converges, it will converge to  a sub-action for $A$.

\bigskip


\begin{lemma} For any $\beta$, we have
$-\|A\| < \frac{1}{\beta}\log\lambda_{\beta} < \|A\|$.

\end{lemma}

 {\bf Proof: }{\bf
 Fix $\beta>0$. Let $\bar x$ be the maximum of $\psi_{\beta A}$ in $\ber$ and $\tilde x$ be the minimum of $\psi_{\beta A}$ in $\ber$. If $\|A\|$ is the uniform norm of $A$, we have
 $$\lambda_{\beta}=\frac{1}{\psi_{\beta A}(\bar x)} \int_M e^{\beta A(a\,\bar x)} \psi_{\beta A}(a\, \bar x)d \nu(a)\leq
 \int_M e^{\beta A(a\,\bar x)} d \nu(a)
 \leq e^{\beta \|A\|} \;,$$
and
 $$\lambda_{\beta}=\frac{1}{\psi_{\beta A}(\tilde x)} \int_M e^{\beta A(a\,\tilde x)} \psi_{\beta A}(a \,\tilde x)d  \nu(a) \geq
  \int_M e^{\beta A(a\,\bar x)}d  \nu(a)
 \geq e^{-\beta \|A\|} \; .$$}
\cqd



The next result can be seen like a kind of measure theoretical version  of the Laplace's method.

\begin{lemma}\label{laplace}
Suppose $W_t:M \to \re$ converges uniformly to $W:M\to \re$, when $t\to +\infty$.
Then \[ \lim_{t \to \infty} \frac{1}{t} \log \int_M e^{t W_t(a)} d\nu(a) = \max_{a\in M}\, W(a).  \]
\end{lemma}

\noindent
{\bf Proof:}
Let $m = \max \{W(a) : a \in M\}$. Let $\bar a \in M$, such that $W(\bar a) = m$.
Given $\varepsilon >0 $, there exist $t_0$ and $\delta$, such that $W_t(a) > m - \varepsilon$, for any $a \in B(\bar a,\delta) \equiv \{a : d(a,\bar a)<\delta \}$ and $t > t_0$.

Therefore, if $t > t_0$, we have  that

$$\int_M e^{t W_t(a)} d\nu(a) \geq \int_{B(\bar a,\delta)} e^{t W_t(a)} d\nu(a) > \nu\big(B(\bar a,\delta)\big) e^{t(m-\varepsilon)},$$
 thus, if $t > t_0$,
\[  \frac{1}{t} \log \int_M e^{t W_t(a)} d\nu(a) >\frac{1}{t}\log\big(\nu\big(B(\bar a,\delta)\big)\big) + m-\varepsilon.\]
Hence
\[  \liminf_{t \to +\infty} \frac{1}{t} \log \int_M e^{t W_t(a)} d\nu(a) \geq m.  \]

The other inequality is analogous, using the fact that: given $\varepsilon$, there exists $t_0$ such that, $W_t(a) < m + \varepsilon$, for any $t>t_0$ and $a\in M$.
 \cqd

\bigskip

\begin{proposition}\label{existsubac}
Given a potential $A$ Hölder continuous, we have\newline
i) $$\lim_{\beta \to \infty}  \frac{1}{\beta}\log\lambda_{\beta}= m(A).$$
ii) Any limit, in the uniform topology,
 $$V:=\lim_{n\to\infty}\frac{1}{\beta_n}\log(\psi_{\beta_n A}),$$
 is  a calibrated subaction for $A$.
\end{proposition}

\noindent{\bf Proof:} Let $\beta_n$ be a subsequence such that the following limit exists: $ \frac{1}{\beta_n}\log\lambda_{\beta_n}\to k$, when $n \to \infty$.
 By taking a subsequence of $\beta_n$ we can assume that  also  there exists $V$ Holder, such that
$\displaystyle V:=\lim_{n\to\infty}
\frac{1}{\beta_n}\log(\psi_{\beta_n A}).$

Given $x \in \cal{B}$, consider the equation
 $$\lambda_{\beta_n}=\frac{1}{\psi_{\beta_n A}(x)} \int_M e^{\beta_n A(a\, x)} \psi_{\beta_n A}(a\,x)d\nu( a) .$$
It follows from lemma \ref{laplace} that, when $n\to \infty$,
$$k= \max_{a \in M} \{A(ax)+  V (ax) - V(x) \} .$$
%

First we show  that $k\geq m(A):$\newline
From the above it follows  that
$$ -V(\sigma(x)) + V(x) + A(x)\leq k.$$
Let $\mu$ be a $\sigma$-invariant  probability measure, then
$$\int_{\ber} A(x) d \mu(x)=\int_{\ber}  \big[-V(\sigma(x)) + V(x) + A(x)\big]\,d \mu(x)\leq k.$$
This implies $m(A)\leq k.$

Now we show that $m(A) \geq k$ :\newline
For any $x\in\ber$ there exist $y=a_x\,x$ such that $\sigma(y)=x$, and
$$ -V(\sigma(y)) + V(y) + A(y)= k.$$
Therefore, the compact set $K= \{y \,:\,-V(\sigma(y)) + V(y) + A(y)= k\}$ is such that, $K' = \cap_n \, \sigma^{-n} (K)$ is non-empty, compact and $\sigma$-invariant. If we consider a $\sigma$-invariant probability measure $\mu$ with support on $K'$, we have that $\int_{\ber} A(y) d \mu(y)=k$. From this follows that $m(A) \geq k$.

From the above arguments $k=m(A)$ is the unique possible accumulation point of the bounded function $\beta \to \frac{1}{\beta}\log\lambda_{\beta\, A}$, then
$$\lim_{\beta \to \infty}  \frac{1}{\beta}\log\lambda_{\beta\, A}= m(A).$$

Moreover, from the above expressions, we can say that any limit of convergent subsequence $\displaystyle\lim_{n\to\infty}
\frac{1}{\beta_n}\log(\psi_{\beta_n A})$ is a calibrated subaction.
\cqd

\bigskip
Now we return to study the Gibbs measures at zero temperature.
 In the case $\mu_{\beta A} \rightharpoonup \mu_{\infty}$, when $\beta \to \infty$ (not just a subsequence), as  we said before,  we have selection of probability  at temperature zero (see \cite{Lep}, \cite{LM}, \cite{LMST} for  general positive results and \cite{CGU} \cite{CH} \cite{van} for negative results).
The  next result uses the variational principle proved in the previous  section and the property that the entropy of an invariant probability  is not positive.

\begin{theorem} Consider a  Hölder potential $A$.
Suppose that for some subsequence
we have $\mu_{\beta_n A} \rightharpoonup \mu_{\infty}$. Then $\mu_{\infty}$
 is a maximizing probability, i.e.,
$$\int _{\ber}A(x)d\mu_{\infty}(x)=m(A) .$$

In the case the maximizing probability for $A$ is unique,  we have selection  of Gibbs probability  at temperature zero.

\end{theorem}

\noindent{\bf Proof.} By definition, $\mu_{\beta_n A} \rightharpoonup \mu_{\infty}$, if and only if,
$$\lim_{n\to\infty} \int_{\ber} w d \mu_{\beta_n A}= \int_{\ber} wd \mu_{\infty}, \,\,\,\, \forall \,w\in \cal{C} . $$
Now using Theorem \ref{pvariacional} and the fact that $h(\mu)\leq 0$, we obtain
\begin{eqnarray*}
  m(A) &=& \lim_{\beta\to\infty}\frac{\log \lambda_{\beta A}}{\beta}=\lim_{n\to\infty} \bigg(\int_{\ber} Ad\mu_{\beta_n A}+\frac{1}{\beta_n}h( \mu_{\beta_n A})\bigg) \\
   &\leq& \lim_{n\to\infty} \int_{\ber} Ad\mu_{\beta_n A}=\int_{\ber} Ad\mu_{\infty} \end{eqnarray*}
Hence, $m(A)\leq \int_{\ber} Ad\mu_{\infty}$. Also, as $\mu_{\infty}$ is a $\sigma-$invariant measure, we have that $m(A)\geq \int_{\ber} Ad\mu_{\infty}$. This implies that $m(A)= \int_{\ber} Ad\mu_{\infty}$.

\cqd

\bigskip

Questions related to the Large Deviation property on the $XY$ model, when $\beta \to \infty$, are considered in \cite{LM}. The existence of a calibrated subaction plays an important role in this kind of result.
\bigskip

We consider now a different kind of question.
From an easy adaptation of Theorem 2.1 in \cite{CGui} one can show:

\begin{proposition}

$$\lim_{n \to \infty} \frac{1}{n} \, \sup_x \, \sum_{j=0}^{n-1} A(\sigma^j(x) )=m(A).$$

\end{proposition}

We consider from now on a potential $A$ which depends on two coordinates $A:S^1 \times S^1 \to \mathbb{R}$ of class $C^\infty$.
A smooth real-valued function on a manifold M is a Morse function if it has no degenerate critical points.

In \cite{AFMT} it is shown that in the $C^\infty$ topology it is  generic  the set of  potentials $A$ such that for any $n$ the function
$\sum_{j=0}^{n-1} A\circ\sigma^j:(S^1)^{n+1}\to \mathbb{R}$ is a Morse function. In this case for each $n$ there exist a finite number of points where the values of
$\frac{1}{n}\,\sum_{j=0}^{n-1} A\circ\sigma^j$ are maximal. Moreover, there exists a positive number $D$ such that for all $n$ the number of critical points of $\sum_{j=0}^{n-1} A\circ\sigma^j$
is smaller than $D^n$.

One can consider for the function $\sum_{j=0}^{n-1} A\circ\sigma^j:(S^1)^{n+1}\to \mathbb{R}$ periodic boundary conditions on $(S^1)^n$.
By adapting the  proof of the above result we get that the   Morse property, in this case, is also true. It follows from the above proposition  that, generically on the potential $A$, the maximizing probability can be approximated by probabilities with support on periodic orbits (which are isolated in $(S^1)^{n+1}).$

Other references related to the topic are \cite{Be} \cite{BG} \cite{FM}.
In a future work we will analyze questions related to zeta functions for generic potentials $A$ (see \cite{PP} and \cite{LM3}).
\bigskip

\section{An application to the non-compact case} \label{non}

An interesting example of application of the above theory is the following: consider  $M_0 = \{z_i, i \in \mathbb{N} \}$ an increasing infinite sequence of points in $[0,1)$
and suppose that $z_\infty:=1=\lim_{i\to\infty} z_i $.
We will also suppose $z_1=0$.
Therefore, each point of $M_0$ is isolated, and there is only one accumulating point $z_{\infty}=1$. 
We consider  the induced euclidean metric then $M=M_0 \cup \{ 1\}$
is a compact set.
The state space $M_0$  can be identified with $\mathbb{N}$, and $M$ has a special point $z_{\infty}=1$ playing the role of the infinity.
Let $\ber_0=M_{0}^{\mathbb{N}}$ and $\ber=M^{\mathbb{N}}$. Note that $\ber_0$ is not compact.

Some results in Thermodynamic Formalism for the shift with countable symbols (see \cite{Sa} \cite{Da}) can be recovered from our previous results  as we will see. We will also study the limit, when the temperature goes to zero, of Gibbs states,  and some results in Ergodic Optimization will be obtained. In particular we will show the existence of sub-actions, under some suitable  hypothesis.  We will get this last result via limit at temperature zero of eigenfunctions
at positive temperature. Results in Ergodic Optimization for this setting appear in
\cite{JMU}, \cite{BG}, \cite{BS}, \cite{Kem}, \cite{Iom}, \cite{Mor}.

\bigskip
\begin{center}
\textbf{Thermodynamic Formalism}
\end{center}

\begin{lemma}\label{extended}
Suppose that $A:\ber_0 \to\mathbb{R}$ is a Hölder continuous potential. Then it can be extended as a Hölder continuous function $A:\ber \to \mathbb{R}$.
\end{lemma}

\noindent
{\bf Proof:}
The extension is a consequence of the fact that any uniformly continuous function can be extended as a uniformly continuous function to the closure of its domain.  It is easy to see that this extension is also Holder continuous.
\cqd

\bigskip
Now let us fix an a-priori measure $\nu:=\sum_{i\in \mathbb{N}}p_i\delta_{z_i}$ on $M$ (or $M_0$), where $p_i>0$ and $\sum_{i\in \mathbb{N}}p_i = 1$.
In fact,
we have that $z_{\infty}=1$ belongs to the support of $\mu$, but is not an atom of $\mu$. All other points of $M$ (i.e. the points of $M_0$) are atoms for $\nu$. 
  On this way for each Hölder continuous potential $A:\ber_{0}\to\R$ we can consider the following Transfer Operator on ${\cal C} (\ber _0)$:
\[\LL_A(w)(x):= \int_M e^{A(ax)}w(ax) d\nu(a) = \sum_{i\in \mathbb{N}}e^{A(z_ix)}w(z_ix)p_i.\]

\begin{proposition}\label{RPF berz} Let $A:\ber_0\to\R$ be a Hölder potential. Then

 (a) there exists a positive number $\lambda_A$ and a positive Hölder function $\psi_A:\ber_0 \to\R$, such that, $\LL_A \psi_A = \lambda_A \psi_A.$\newline
   If we consider
the normalized potential
$\bar A = A+\log  \evc-\log \evc \circ \sigma -\log \evl$, then

(b)  there exists an unique fixed point $\mu_{A}$ for $ {\cal L }_{\bar A}^*$,  which  is a $\sigma$-invariant probability measure on $\ber_0$.

(c) the measure $$\rho_A= \frac{1}{\psi_A} \,\, \mu_{A}$$ satisfies ${\cal  L}_A^* (\rho_A)= \lambda_A \rho_A$. Therefore, $\,\rho_A$  is an eigen-measure for ${\cal  L}_A^*$.

(d) for any Hölder function $w:{\ber_0} \to \mathbb{R}$, we have that, in the uniform convergence topology,
 $$\frac{{\cal L}_A^n (w)}{(\lambda_A)^n} \to \, \psi_A \int_{\ber_0}\,w \,d \rho_A,$$
 and $${\cal L }_{\bar A}^n\omega \rar \int_{\BB_0} \omega d\mu_{A}\,.$$
\end{proposition}

\medskip
\noindent
{\bf Proof:}
Using  lemma \ref{extended} we can extend $\LL_A$ to $\cal{C}(\ber)$. From Theorem \ref{RPF_eig_nonnorm} we obtain $\lambda_A>0$ and $\psi_A>0$ Hölder continuous such that
\[   \sum_{i}e^{A(z_ix)}\psi_A(z_ix)p_i = \lambda_A.\psi_A(x) \ \ \ \ \forall x\in \ber.\]
In particular,  the restriction $\psi_A:\ber_0\to\R$ will satisfy also the expression
\[  \sum_{i}e^{A(z_ix)}\psi_A(z_ix)p_i = \lambda_A.\psi_A(x) \ \ \ \ \forall x\in \ber_0.\] This proves item (a).\newline
Consider $\bar A = A+\log  \evc-\log \evc \circ \sigma -\log \evl$.

In order to prove  item (b) we observe that from Theorem \ref{RPF_normalized} there exists $\mu_A$ on $\ber$ satisfying  item (b). We want prove that $\mu_A(\ber -\ber_0)=0$, or equivalently $\mu_A(\ber_0)=1$. On this way we only need to show that $\mu_A(\{x \in \ber \,: \, x_1=1\})=0$,
because for all $n\geq 1$, we have  $$\mu_A(\{x \in \ber \,: \, x_n=1\})=
\mu_A(\sigma^{-n+1}(\{x \in \ber \,: \, x_1=1\})),$$  therefore, we will have
$$ \mu_A(\ber - \ber_0) \leq \sum_{n=1}^{\infty} \mu_A(\{x \in \ber \,: \, x_n=1\})=0.$$
To prove that $\mu_A(\{x \in \ber \,: \, x_1=1\})=0$, we fix  $\varepsilon>0$, and we consider  a Hölder function $w_{\varepsilon} $, such that,
$\chi_{\{x\in \ber\,:\, x_1=1\}} \leq w_{\varepsilon} \leq 1$
and $w_{\varepsilon}(x)=0$, if $x_1 < 1-\varepsilon$.
Then, using Theorem \ref{RPF_normalized} item c)
\begin{eqnarray*}
   & & \mu_A(\{x \in \ber \,: \, x_1=1\}) \leq \int_{\ber} w_{\varepsilon} d\mu_A = \lim_{n\to+\infty} \LL_{\bar A}^n w_{\varepsilon} (0^{\infty}) \\
   & & = \lim_{n\to+\infty}\sum_{i_1:z_{i_1} >1-\varepsilon} p_{i_1} \sum_{i_2,...,i_n} e^{S_{\bar A}^n(z_{i_1}z_{ i_2} ...z_{ i_n} 0^{\infty})} w_{\varepsilon}(z_{i_1} z_{i_2} ... z_{i_n} 0^{\infty}) p_{i_2}...p_{i_n} \\
  & & \leq \lim_{n\to+\infty} \sum_{i_1: z_{i_1} >1-\varepsilon} p_{i_1} e^{\|\bar A\|} \sum_{i_2,...,i_n} e^{S_{\bar A}^{n-1}(z_{i_2} ... z_{i_n} 0^{\infty})} p_{i_2}...p_{i_n} =
\sum_{i_1 : z_{i_1}>1-\varepsilon} p_{i_1} e^{\|\bar A\|},
\end{eqnarray*}
where we used in the last equation the normalization property.

Now the claim follows easily when we use the fact that the a-priori measure is supported on $[0,1)$, which makes $$\sum_{i_1 : z_{i_1}>1-\varepsilon} p_{i_1} \to 0 \ \ \ \ \text{when} \ \ \ \ \varepsilon \to 0.$$

The items (c) and (d) follow when we restrict to $\ber_0$ the result of Theorem \ref{RPF_normalized}.
\cqd

\bigskip

Now let us compare this setting with some results contained in \cite{Sa}. The operator $\LL_A$ can be written as
\[\LL_A(w)(x)= \sum_{i}e^{A(z_ix)}w(z_ix)p_i = \sum_{i}e^{A(z_ix)+\log(p_i)}w(z_ix),\]
that is, the Classical Ruelle Operator with potential $B:=A+\log(P)$, where $P(y_1,y_2,y_3,...)=P(y_1)=p_{i}$, if, $y_1=z_i$. We denote this operator by $L_B$, or, $L_{A+\log(P)}$.

Clearly $(A+\log(P))(z_i,y_2,y_3,...)\to-\infty$, when $i\to+\infty$, because $p_i \to 0$, when, $i \to +\infty$. Furthermore, if we define  $$Var_n (B) = \sup\{|B(x)-B(y)|:x_1=y_1,...,x_n=y_n\},$$
then, there exists $C>0$, such that, $Var_n( B) \leq C\frac{1}{2^{n\alpha}}$, for any $n\geq 1$. This means that $B$ is \textbf{locally Hölder continuous} (see \cite{Sa}).

Define
\[Z_{n}(B,a) := \sum_{\scriptsize{ \begin{array}{c}\sigma^{n}(y)=y\\y_1=a\end{array}}\normalsize} e^{S_nB(y)} .\]

\begin{proposition}
Fix $a\in M_0$, then, there exists a constant $M_a$
and an integer $N_a$, such that, for any $n > N_a$:
\[\frac{Z_{n}(B,a)}{(\lambda_A)^{n}} \in [M_a^{-1},M_a]\]
\end{proposition}

\noindent
{\bf Proof:}
Let $x=a^{\infty}=(a,a,a,a,...)$. When we apply item d) of Proposition  \ref{RPF berz}   for $w\equiv 1$ we get
\[\frac{{\cal L}_A^n (w)(x)}{(\lambda_A)^n} \to C>0. \]
Remember that we denote by $a^n=(a_n,...,a_1)\in M_0^n$, then there exist $M_1>0$ and $N_a>0$ such that for $n\geq N_a$:
\begin{eqnarray*}
   & & \frac{1}{(\lambda_A)^{n+1}}\sum_{a^n}e^{S_nB(a^nx)} \\
   & & =\frac{1}{(\lambda_A)^{n+1}}\sum_{a^n}e^{S_nA(a^nx)+\log(P(a_n)...P(a_1))} \in [M_1^{-1},M_1].
\end{eqnarray*}
  Let $y$ be the periodic point with period $n+1$ obtained by the successive concatenation of $(a,a_n,...,a_1)$ and let $z=\sigma(y)$.
We have
\begin{eqnarray*}
  |S_nB(a^nx)-S_nB(z)| &=& |S_nB(a_n,...,a_1,a,a,a,...)-S_nB(a_n,...,a_1,a,a_n,a_{n-1},...)| \\
   &=& |S_nA(a_n,...,a_1,a,a,a,...)-S_nA(a_n,...,a_1,a,a_n,a_{n-1},...)| \\
   &\leq& C\left(\frac{1}{2^{\alpha}}+\frac{1}{2^{2\alpha}}+...+\frac{1}{2^{n\alpha}}\right).d(x,y)^{\alpha}.
\end{eqnarray*}

Using the fact  that $\ber_{0}$ has finite diameter we obtain a constant $M_2>0$, such that,
\[|S_nB(a^nx)-S_nB(z)| \leq M_2, \ \ \ \ \ \ \forall \,n\in\mathbb{N},\, a^n\in M_0^n.\]
Furthermore, using the property  $\sigma^{n}(z) = y $, we conclude that $B(\sigma^{n}z)=A(\sigma^{n}(z))+\log(P(a))$ is bounded independently of $z$. Therefore, there exists a constant $M_3>0$, such that,
\[|S_nB(a^nx)-S_{(n+1)}B(z)| \leq M_3, \ \ \ \ \ \ \forall \,n\in\mathbb{N},\, a^n\in M_0^n.\]

Note  that $S_{(n+1)}B(z)=S_{(n+1)}B(y)$. Then,
\[\sum_{a^n}e^{S_nB(a^nx)-M_3} \leq \sum_{\scriptsize\begin{array}{c}\sigma^{n+1}(y)=y\\ y_1=a\end{array}}\normalsize e^{S_{(n+1)}B(y)} \leq
\sum_{a^n}e^{S_nB(a^nx)+M_3}\]
Therefore, when $n\geq N_a$, we get
\[\left(\frac{1}{(\lambda_A)^{n+1}}\sum_{\scriptsize\begin{array}{c}\sigma^{n+1}(y)=y\\ y_1=a\end{array}}\normalsize e^{S_{(n+1)}B(y)}\right) \in [(M_1e^{M_3})^{-1},(M_1e^{M_3})].\]
Choosing $M_a=M_1e^{M_3}$ we conclude the proof.
\cqd

\bigskip

 In this way, we can say that $B=A+\log(P)$ is \textbf{positive recurrent} (see \cite{Sa} Definition 2). Following \cite{Sa} Theorem 4 we get a Ruelle-Perron-Frobenius Theorem (as in Theorem \ref{RPF berz} above). It follows from the above proposition that $\lambda_A$ is the Gurevic pressure of $B$ (see \cite{Sa} definition 1).

We would like to  point out some differences on the topology considered in our setting with the classical one used in the theory of Thermodynamic  Formalism with state space $\mathbb{N}$. The set $M_0^{\mathbb{N}}$ can be identified with  $\mathbb{N}^{\mathbb{N}}$, but the metric space $M_0^{\mathbb{N}}$ is different from the metric space $\mathbb{N}^{\mathbb{N}}$ with the discrete product topology. Here, we consider a distance (induced in the subset $M_0 \cup \{ z_{\infty}\} \subset [0,1]$), such that, for any two points $x=(x_1,x_2,...),\,y=(y_1,y_2,...) \in M_0^{\mathbb{N}}$
\[d(x,y) = \sum_{n\in\mathbb{N}}\frac{1}{2^{n}}d_{[0,1]}(x_n,y_n).\]
 On the other hand,  the metric considered in \cite{Sa} is of the form: for  two points $x,y \in \mathbb{N}^{\mathbb{N}}$
\[\widetilde{d}(x,y) = \frac{1}{2^{n}}, \,\,\,\, \text{if}\,\,\,\, x_1=y_1,...,x_{n-1}=y_{n-1},\,\,x_{n}\neq y_{n}.\]

Using that the
 diameter of $[0,1]$ is one,   it follows that $d(x,y) \leq \widetilde{d}(x,y)$.
In particular,  any convergent sequence on the metric $\widetilde{d}$ is a convergent sequence on the metric $d$, and any continuous/Hölder function $A$ for the metric $d$ is a continuous/Hölder function for the metric $\widetilde{d}$.
But the same is not true in the opposite direction. This is a subtle question. Results in  \cite{Sa} and here  are
obtained under slight  different hypothesis.  Anyway, in physical applications  this  is probably    a not very important point.

Considering the dual space, it follows from the relation $d(x,y) \leq \widetilde{d}(x,y)$ that any open set for the metric $d$ is an open set for the metric $\widetilde{d}$. Then, the Borel sigma-algebra generated by $d$ is contained in the Borel sigma-algebra generated by $\widetilde{d}$. on the order hand, the cylinder sets \cite{Sa} are closed sets for the metric $d$, therefore, they belong to the sigma-algebra generated by $d$. In this way, the Borel sigma-algebra generated by $d$, or, by $\widetilde{d}$, is the same.

\begin{center}
\textbf{Ergodic Optimization}
\end{center}

Once more we point out that the concepts of sub-action and maximizing measures do not involve an a-priori measure. On this way the statement of the next theorem does not use any a-priori probability. On the other hand, some condition on $A$ must be assumed in order  to obtain positive results in Ergodic Optimization. For example: the potential $A:\ber_0\to\R$, given by $$A(x)=-\sum_{i=1}^{\infty}\frac{1}{2^{i}}d(x_i,1),$$
does not have maximizing measures on $\ber_0$.

We need some hypothesis on $A$ in such way it prevents ``the mass to go to infinity". Under some appropriate and natural conditions  on $A$, we will obtain below  the existence of calibrated sub-actions and maximizing measures. In the proof, we consider a limit involving a subsequence of eigenfunctions of $\LL_{\beta A}$, when the temperature $\frac{1}{\beta}$ goes to zero.

\begin{theorem}Suppose that $A:\ber_0 \to\mathbb{R}$ is a Hölder continuous potential. Consider the Hölder continuous extension $A:\ber \to \mathbb{R}$.
If the extension satisfies:
\begin{equation}\label{hipotese A} A(x_1,...,x_{n-1},1,x_{n+1},x_{n+2},...) < A(x_1,...,x_{n-1},0,x_{n+1},x_{n+2},...) \end{equation}
for any $n\in\mathbb{N}$ and $x_i\in M$, then:\newline
a) A has a calibrated subaction $V$ on $\ber_0$, that means: for any $x\in \ber_0$,
\[0=\max_{a\in M_0}\left(A(ax)+V(ax)-V(x)-m(A)\right).\]
b) Any maximizing measure for $A$ has support on $\ber_0$.
\end{theorem}

\noindent
{\bf Proof:}
We define an a-priori probability measure on $M_0$ by the expression
$$\nu = \sum_{i \in \mathbb{N}} p_i. \delta_{z_i}  ,$$ where $\displaystyle\sum_{i\in \N } p_i =1 $ and $p_i > 0$.

{\it Claim:} Denote by  $\psi_{\beta A}$  the maximal eigenfunction for the Ruelle operator given by the potential $\beta A$ and the a-priori measure $\nu$. Let $V$ be the limit of some subsequence $\frac{1}{\beta_n} \log \psi_{\beta_n A}$. Then, we have
 $$ (i) \;\;\;\;    \psi_{\beta A}(1,x_1,x_2,x_3,...)<\psi_{\beta A}(0,x_1,x_2,x_3,...) \;\;\;\forall x\,\in \ber,$$
 $$(ii)\;\;\;\; V(1,x_1,x_2,x_3,...)\leq V(0,x_1,x_2,x_3,...) \;\;\;\;\forall x\,\in \ber.$$
\newline Proof of (i):
 We know by item (c) of the Theorem \ref{RPF_normalized} that
 $$\frac{{\cal L}_{\beta A}^n (1)(x)}{(\lambda_{\beta A})^n} \to \, \psi_{\beta A}(x) \int_{\ber}\,1\,d \rho_{\beta A}, \;\;\;\;\forall x\in\ber.$$
Hence,
$$ \rho_{\beta A}(\ber)\psi_{\beta A}(x) =\lim_{n\to\infty} \int_{M^n}  e^{S_n{\beta A}(a^n x)-n\log \lambda_{\beta A}}  d\nu^n(a^n),\;\;\;\;\forall x\in\ber.$$
Suppose now that $z=(1,x_1,x_2,x_3,...)$ and $y=(0,x_1,x_2,x_3,...)$, with $x\in \ber$,  
then, using \eqref{hipotese A}, we get
$$\int_{M^n}  e^{S_n{\beta A}(a^n z)-n\log \lambda_{\beta A}}  d\nu^n(a^n)< \int_{M^n}  e^{S_n{\beta A}(a^n y)-n\log \lambda_{\beta A}}  d\nu^n(a^n),$$ which implies
$\psi_{\beta A}(z)<\psi_{\beta A}(y)$.\newline
Proof of item (ii): By hypothesis
 $V$ satisfies $$V(x)=\lim_{n\to\infty}
\frac{1}{\beta_n}\log(\psi_{\beta_n A}(x)),$$
 as  the $\log$ function is monotone, we get  $V(z)\leq V(y)$,
which finishes the proof of the claim.

\medskip

Now 
let $R_{-}= A+V-V\circ\sigma - m(A)$, where $V$ was defined above. We know, by the sub-action equation, that $R_{-}\leq 0$. From ($\ref{hipotese A}$) and the above claim, we get that for any $x\in \ber_0$:
\begin{eqnarray*}
 R_{-}(1x) &=& A(1x)+ V(1x)- V(x)-m(A) \\
   &<& A(0x)+ V(0x)- V(x)-m(A) = R_{-}(0x)\leq 0.
\end{eqnarray*}
Using now the fact  that $V$ is a calibrated subaction (on $\ber$) we conclude the proof of $(a)$, because last inequality shows that
$$\max_{a\in M_0}\left(A(ax)+V(ax)-V(x)-m(A)\right)=
\max_{a\in M}\left(A(ax)+V(ax)-V(x)-m(A)\right).$$

We point out that the fact that the extension has a maximizing measure is a consequence of the compactness of $\ber$. In order to prove (b) we will fix a  $x\in \ber \verb"\" \ber_0$ and prove that $x$ does not belong to the support of any maximizing measure.

Note that if $R_{-}(\sigma^{k-1}(x))<0$, then $x$ does not belong to the support of the maximizing probability $\mu$. Indeed,

$$ \int_{\ber} R_{-} \circ \sigma^{k-1} d\mu = \int_{\ber} R_{-} d\mu = \int_{\ber} A d\mu - m(A) =0,  $$
which, combined to the continuity of $R_{-}\leq0$,  proves that $R_{-}\circ \sigma^{k-1}$ vanishes at the support of $\mu$.
Now we will prove that $R_{-}(\sigma^{k-1}(x))<0$. So, let $k \in \mathbb{N}$ be such that $x_k=1$ and $x_l<1$, $ \forall 1\leq l  <k$.
Let $y \in \ber$ be given by $y_i=x_i$, if $i \neq k$, and $y_k=0$.

We have that $\sigma^{k-1}(x) = (1,x_{k+1},x_{k+2},...)$ and $\sigma^{k-1}(y) = (0,x_{k+1},x_{k+2},...)$, therefore
\begin{eqnarray*}
  R_{-}(\sigma^{k-1}(x)) &=& A(\sigma^{k-1}(x))+ V(\sigma^{k-1}(x))- V(\sigma^{k}(x))-m(A) \\
   &<& A(\sigma^{k-1}(y))+ V(\sigma^{k-1}(y))- V(\sigma^{k}(y))-m(A) = R_{-}(\sigma^{k-1}(y))\leq 0,
\end{eqnarray*}
where we used  above the hypothesis \eqref{hipotese A}, item (ii) of the claim, and also $\sigma^{k}(x)=\sigma^{k}(y)$.
\cqd

\bigskip

In the proof of the above theorem we show the following result:

\begin{corollary}
Given an a-priori probability measure $\nu = \sum_{i=1}^{\infty}p_i\delta_{z_i}$, $p_i>0$, under the hypothesis of the above theorem, then,
there exists a subsequence $\{\beta_n\}$ and a Hölder continuous function $V:\ber_0\to\R$, such that,
\[\frac{1}{\beta_n}\log(\psi_{\beta_nA}) \to V\]
uniformly on $\ber_0$. Furthermore, any possible limit $V$ is a calibrated sub-action for $A$ on $\ber_0$.
\end{corollary}

\medskip

An example of potential satisfying the hypothesis of the above theorem is given by
\[A(x) = -d(x,0^{\infty})\]
where $0^{\infty}=(0,0,0,0...)$. Note that in the claim of  hypothesis $(\ref{hipotese A})$ we can change $0$ by any fixed $z_i$, and the result we get  will be the same.

\bigskip

It is important to remark that when the temperature changes, then,  the operator varies in a different way of what happens in the classical sense: for a fixed $\beta>0$, we have:
\[\int_M e^{\beta\,A(ax)}w(ax) d\nu(a) = \sum_{i}e^{\beta\,A(z_ix)}w(z_ix).p_i =\]
\[ \sum_{i}e^{\beta\,A(z_ix)+\log(P(z_i))}w(z_ix) .\]
Then, (in this work) the main eigenvalue and eigenfunction, respectively,  $\lambda_{\beta}>0$ and $\psi_{\beta}>0$, are associated to $\beta(A+\frac{\log(P)}{\beta})$ (in the setting of the Classical Ruelle Operator).

 In this way we can think of the  function $z_i\to\frac{\log(P(z_i))}{\beta}$ as a perturbation of the potential $A$, that goes to zero when $\beta \to\infty$ (but not uniformly).



\section{The differentiable structure and the involution kernel} \label{dif}

We consider in this section the  $XY$ model. This is the case where $M=S^1$, and the a-priori measure is the Lebesgue measure on the circle.
$(S^1)^\mathbb{N}$ has a differentiable structure.

\textbf{ We know that, in the case where the potential $A$ is Holder, the eigenfunction $\psi_A$ is also Holder and belongs to the same Holder class. The main result of this section is theorem \ref{differentiability}, where we prove that, under mild assumptions concerning the differentiability of $A$ (see definition \ref{def-classofdiff}), the associated eigenfunction $\psi_A$ is differentiable in each coordinate $x_j$ of $x$.}


In this setting, we point out that  in  \cite{LMST} it is analyzed several questions which involve differentiability for potentials which depend  just on two coordinates. Here we consider more general potentials.

Let $\ber^*=\{(...,y_2,y_1)\in (S^1)^{\mathbb{N}}  \}$, and we denote by the pair
$$( y|x) =(...,y_2,y_1|x_1,x_2...)  ,$$
the general  element of $\hat{\ber}:=\ber^*\times\ber=(S^1)^{\mathbb{Z}}$,\textbf{ the natural extension of $\ber$. Here we will follow the ideas of \cite{BLT}.}

We denote by $\hat{\sigma}$ the shift on $\hat{\ber}$,  i.e.
$$\hat{\sigma}(  ...,y_2,y_1|x_1,x_2,...)= ( ...,y_2,y_1,x_1|x_2,x_3,...).$$

\begin{definition}
Let $A:\ber\to \mathbb{R}$ be a continuous  potential (considered as a function on $\hat{\ber}$). A continuous function   $W:\hat{\ber}\to\re$ is called an involution kernel, if
$$A^*:=A\circ\hat{\sigma}^{-1}+W\circ\hat{\sigma}^{-1}-W $$ depends only on the variable $y$.

\end{definition}


 The involution kernel is not unique.
 \medskip

Let us fix $x'\in\ber$ and $A$ a Hölder continuous potential, then we define
\begin{equation}\label{Wkernel} W(y|x)=\sum_{n\geq 1}A(y_n,...,y_1,x_1,x_2,...)-A(y_n,...,y_1,x_1',x_2',...).  \end{equation}
 An easy calculation shows that $W(y|x)$ is a involution kernel (see \cite{BLT}).

\textbf{The Ruelle-Perron operator $\mathcal{L}_A$ gives two important informations: the eingenmeasure $\rho_A$ and the eingenfunction $\psi_A$. As in \cite{BLT}, we can use the involution Kernel in order to obtain $\psi_A$, if we know the eigenmeasure of the Ruelle-Perron operator associated to $A^*(y)$ (see proposition \ref{psi W kernel}). We
 will show that the involution kernel allows one  to differentiate $\psi_A$ with respect to each coordinate $x_j$ of $x$, using the expression of $\psi_A$ given in proposition \ref{psi W kernel}
 (see theorem \ref{differentiability}).}
 \medskip

Let $A:\ber\to \mathbb{R}$ be a Hölder continuous  potential  and  $W:\hat{\ber}\to\re$  an involution kernel, then
for any $a\in S^1$, $x\in\ber$ and
$y\in\ber^*$, we have
\begin{equation}\label{eqWkernel}
(A^*+W)( ya|x)=(A+W)( y|ax).
\end{equation}
Questions related to Ergodic Transport Theory and the involution kernel are analyzed in \cite{LOT}, \cite{CLO} and  \cite{LM2}.

Remember that
$$m(A) = \sup_{\mu\,\text{is}\,\sigma-\text{invariant } \, } \,\int_{\ber} A d \mu ,$$
and define
$$m(A^*)= \sup_{\mu\,\text{is}\,\sigma^*-\text{invariant } \, } \,\int_{\ber^*} A^* d \mu.$$
\medskip

The next result is an adaptation to the present setting of a result in \cite{BLT}.

\begin{lemma}\label{label9}

Let $\mathcal{L}_A$ and $\mathcal{L}_{A^*}$ be the Ruelle operators defined on
$\ber$ and $\ber^*$, and $W(y|x)$ an involution kernel.

Then, for any $x\in\ber$,
$y\in\ber^*$, and any function
$f:\hat{\ber}\to\mathbb{R}$
\begin{equation}\label{eq L}
\mathcal{L}_{A^*}\Big(f(\cdot|x)
\, e^{W(\cdot|x)}\Big)(y)
=\mathcal{L}_{A}\Big(f\circ\hat\sigma(y|\cdot)
\, e^{ W(y|\cdot)}\Big)(x).
\end{equation}

\end{lemma}

{\bf Proof:}
Under our notation we write $A(y|x)=A(x)$ and $A^*(y|x)=A^*(y)$. Consider $x\in\ber$,
$y\in\ber^*$ fixed, then by the definition of $\mathcal{L}_A$ and $\mathcal{L}_{A^*}$, and equation \eqref{eqWkernel}, we obtain

\noindent
\begin{align*}
\mathcal{L}_{A^*}\big(f(\cdot|x)\,
\, e^{ W(\cdot|x)}\big)(y)&=\int_{S^{1}} f( ya|x)\, e^{
\big(A^*( ya)+W( ya|x)\big)}\,d\nu(a)\\
&=\int_{S^{1}} \,f\circ\hat\sigma(y|ax)
\, e^{\big(
A(ax)+W(y|ax)\big)}\,d\nu(a) \\
&=\mathcal{L}_{A}\big(f\circ\hat\sigma(y|\cdot)
\, e^{ W(y|\cdot)}\big)(x)
\end{align*} \cqd

Let $\rho_A$ and $\rho_{A^*}$ the eigenmeasures for $\mathcal{L}^*_A$ and $\mathcal{L}^*_{A^*}$, given in Theorem \ref{RPF_normalized}.
Suppose $c$ is such that
$\iint_{\ber\times\ber^*} \, e^{ W(y|x)-c}\, d\rho_{A^*}(y) d\rho_A(x)=1$.


\begin{proposition}\label{psi W kernel}

Suppose $K(y|x)=e^{ W(y|x)-c}$. Then,
$$d\,\hat{\mu}_A= K(y|x)\,d\rho_{A^*}\,(y) d\rho_{A}(x)$$
is invariant for $\hat{\sigma}$ and is the natural extension of the Gibbs measure $\mu_A$.

The function $\psi_A(x)=\int_{\ber^*}\!K(y|x)\,d\rho_{A^*}(y)$ is  the main eigenfunction for $\mathcal{L}_{A},$ and
 the function $\psi_{A^*}(y)=\int_\ber\!K(y|x)\,d\rho_{A}(x)$ is  the main eigenfunction for $\mathcal{L}_{A^*}. $ Furthermore $\lambda_A=\lambda_{A^*}.$

\end{proposition}

{\bf Proof:}
We denote by
$K(y|x)=e^{ W(y|x)-c},$ and we define a positive function $\psi$, by the expression $\psi(x)=\int_{\ber^*}\!K(y|x)\,d\rho_{A^*}(y)$. In order to prove that $\psi$ is an eigenfunction for $\mathcal{L}_{A},$ we remember  that $\mathcal{L}^*_{A^*}(\rho_{A^*})=\lambda_{A^*} \rho_{A^*}$, hence
\begin{align*}
\psi(x)&=\int_{\ber^*}\!K(y|x)  \, d\Big(\frac{1}{\lambda_{A^*}}\mathcal{L}^*_{A}(\rho_{A^*})\Big)(y)=\int_{\ber^*} \frac{1}{\lambda_{A^*}}\mathcal{L}_{A^*}\big(
K(\cdot|x)
\big)(y)\,d\rho_{A^*}(y)\\
&=\int_{\ber^*} \frac{1}{\lambda_{A^*}}\mathcal{L}_{A}\big(
K(y|\cdot)
\big)(x)\,d\rho_{A^*}(y)=\frac{1}{\lambda_{A^*}}\mathcal{L}_{A}(\psi)(x),
\end{align*}
where in the third equality we have used equation \eqref{eq L} with $f=1$. This means that $\psi$ is a positive eigenfunction for  $\mathcal{L}_{A},$ now using Proposition \ref{proponly} we get that $\psi=\psi_A$ and $\lambda_{A^*}=\lambda_A$. The proof for the case of $\psi_{A}^*$  is similar.

By the same arguments used above, for any bounded Borel $f:\ber^* \times \ber \to \mathbb{R}$, we have
\begin{eqnarray*}
   & &\int_{\ber} \int_{\ber^*} f\circ\hat\sigma(y|x)K(y|x)
d\rho_{A^*}(y) d\rho_A(x) \\
   & & =\int_{\ber^*}\!d\rho_{A^*}(y)\int_{\ber} \frac{1}{\lambda_{A}}\mathcal{L}_{A}\big(
f\circ\hat\sigma(y|\cdot)K(y|\cdot)
\big)(x)\,d\rho_A(x) \\
   & & =\int_{\ber}\!d\rho_A(x)\int_{\ber^*} \frac{1}{\lambda_{A^*}}\mathcal{L}_{A^*}\big(
f(\cdot|x)K(\cdot|x)
\big)\,d\rho_{A^*}(y) \\
   & & =\int_{\ber} \int_{\ber^*} f(y|x)\, K(y|x) d\rho_{A^*}(y) d\rho_A(x),
\end{eqnarray*}
hence $d\hat{\mu}_A= K(y|x)\,d\rho_{A^*}\,(y) d\rho_{A}(x)$
is invariant for $\hat{\sigma}$.

Finally, let us prove that $\hat{\mu}_A$ is the natural extension of $\mu_A$.
Given a function $f(x)$ we get that
\begin{eqnarray*}
  \int_{\ber} \int_{\ber^*} f(x) d\hat\mu_A(y,x) &=& \int_{\ber} \int_{\ber^*} f(x)\, K(y|x)\,d\rho_{A^*}(y) d\rho_A(x) \\
   &=& \int_{\ber} f(x) \, \psi_A (x) d \rho_A (x)=\int_{\ber} f(x) d \mu_A (x).
\end{eqnarray*}
Therefore, the measure $d\hat\mu_A(y,x)=K(y|x)\,d\rho_{A^*}(y) d\rho_A(x)$  projects onto $\mu_A$ and $\mu_{A^*}$ (by the same arguments). The probability $\hat\mu_A$ is therefore the natural extension of $\mu_A$.
\cqd

\bigskip

\noindent\textbf{Remark:} Note that, as $\lambda_{\beta A}=\lambda_{\beta A^*}$, we have that $\displaystyle m(A)=\lim_{\beta\to\infty}\frac{1}{\beta}\log\lambda_{\beta A}=m(A^*)$.

\bigskip

{\bf
\begin{definition}\label{def-classofdiff}
Suppose that $A$ is Lipschitz continuous and satisfies the both conditions described below:

(a) $A$ is differentiable in each coordinate $x_j$ of $x\in\ber$,

(b) given $\varepsilon>0$, there exists $H_{\varepsilon}>0$, such that, for all $x$, if $|h|<H_{\varepsilon}$, then
\begin{equation}\label{deriv A}\bigg|\frac{A(x+he_j)-A(x)}{h} -D_jA(x)    \bigg|\leq \frac{\varepsilon}{2^j} \,\,\, \forall j\in\mathbb{N},  \end{equation} where $D_jA(x)$ denote the derivative of $A$ with respect to the $j$-th coordinate.

We will denote the  class of such potentials by $\mathcal{D}$.
\end{definition}


The potential $A(x)= \sum_{n=1}^\infty \frac{1}{2^n} \sin( x_n+1/2^n )$ belongs to the class $\mathcal{D}$ . Moreover, any potential which depends on finite coordinates and is of class $C^2$ belongs to $\mathcal{D}.$


\begin{proposition} Suppose that $A$ belongs to the class $\mathcal{D}$. Given an involution kernel $W$ we have that for any $j$
$$ \frac{\partial }{\partial x_j}W(y|x)   =\sum_{n\geq 1}D_{n+j}A(y_n,...,y_1,x_1,x_2,...)\,.$$
\end{proposition}
}

{\bf Proof:} {\bf Let us first prove that the sum in the right hand side is convergent. Indeed, using that $A$ is Lipschitz,  there exist   $K>0$ such that
$$|A(x)-A(\tilde x)|<K d(x,\tilde x).$$
If we denote  $x^j_h=x+he_j$ and
$y^n x=(y_n,...,y_1,x_1,x_2,...)$, then we have
$$ \bigg|\frac{A(y^n x+h e_{j+n})-A(y^n x)}{h}     \bigg|\leq \frac{K}{h}\frac{d_{S^1}(x_j+h,x_j)}{2^{n+j}}=\frac{K}{2^{n+j}}.$$ Using that $A$ belongs to $\mathcal{D}$, we get that, given $\varepsilon >0$, there exists $H_{\varepsilon}>0$  such that, for each $|h|<H_{\varepsilon}$, we have
$$|D_{n+j}A(y^n x)|\leq \bigg|  \frac{A(y^n x^j_h)-A(y^n x)}{h}     \bigg|+\frac{\varepsilon}{2^{n+j}}\leq \frac{K}{2^{n+j}} + \frac{\varepsilon}{2^{n+j}},  $$ which implies
$$\sum_{n\geq 1}|D_{n+j}A(y^n x)|< \frac{K+\varepsilon}{2^j} <\infty.$$
Now, we will prove the proposition. 
\begin{eqnarray*}
& & \bigg|\frac{W(y|(x+he_j))-W(y|x)}{h}      - \sum_{n\geq 1}D_{n+j}A(y^n x) \bigg|\\
& & =\bigg|\frac{1}{h} \Big(\sum_{n\geq 1} A(y^n x+he_{j+n})-A(y^n x) \Big)- \sum_{n\geq 1}D_{n+j}A(y^n x) \bigg|\\
& & = \bigg|\sum_{n\geq 1}\bigg(\frac{A(y^n x+he_{j+n})-A(y^n x)}{h} - D_{n+j}A(y^n x)  \bigg)\bigg|   \leq \sum_{n\geq 1}\frac{\varepsilon}{2^{n+j}}=\frac{\varepsilon}{2^{j}},
\end{eqnarray*}
for all $|h|\leq H_{\varepsilon}$, as $A$ belongs to $\mathcal D$.

\cqd

From the final part of the last proof we have
that, for all $|h|\leq H_{\varepsilon}$, and for all $x$ and $y$,

\begin{equation}\label{difW}
\bigg|\frac{W(y|(x+he_j))-W(y|x)}{h}  -  \frac{\partial W (y|x)}{\partial x_j}   \bigg|\leq
\frac{\varepsilon}{2^{j}}\,.
\end{equation}


{\bf

\begin{theorem}\label{differentiability}
Let $\psi_A(x)=\int_{\ber^*}\!e^{ W(y|x)-c}\,d\rho_{A^*}(y)$, and suppose $A$ belongs to the class $\mathcal{D}$.  Then,  the eigenfunction $\psi_A$ is differentiable in each coordinate $x_j$. Moreover,
$$\frac{\partial }{\partial x_j}\psi_A(x)=\int_{\ber^*}   e^{ W(y|x)-c} \sum_{n\geq 1}D_{n+j}A(y_n,...,y_1,x_1,x_2,...)\,d\rho_{A^*}(y) . $$

\end{theorem}

}



{\bf Proof:} Consider $j \in \mathbb{N}$. We have
\begin{eqnarray*}
& & \bigg| \frac{\psi_A(x+he_j)-\psi_A(x)}{h}- \int_{\ber^*} e^{W(y|x)-c}\;\;\frac{\partial W(y|x)}{\partial x_j} \;d\rho_{A^*}(y)
\bigg| \\
& & =  \bigg| \int_{\ber^*} \bigg(\frac{e^{W(y|x+he_j)-c}-e^{W(y|x)-c}}{h}- e^{W(y|x)-c}\;\;\frac{\partial W(y|x)}{\partial x_j} \bigg) \;d\rho_{A^*}(y)
\bigg| \\
& & =  \bigg| \int_{\ber^*} e^{W(y|x)-c} \bigg(\frac{e^{W(y|x+he_j)-W(y|x)}-1}{h}- \;\;\frac{\partial W(y|x)}{\partial x_j} \bigg) \;d\rho_{A^*}(y)
\bigg| \\
& &  \leq   \bigg| \int_{\ber^*} e^{W(y|x)-c} \bigg(\frac{e^{W(y|x+he_j)-W(y|x)}-1}{h}- \;\;\frac{W(y|x+he_j)-W(y|x)}{h} \bigg) \;d\rho_{A^*}(y)
\bigg|\\
& &  + \bigg| \int_{\ber^*} e^{W(y|x)-c} \bigg(\frac{W(y|x+he_j)-W(y|x)}{h}- \;\;\frac{\partial W(y|x)}{\partial x_j} \bigg) \;d\rho_{A^*}(y)
\bigg|.
\end{eqnarray*}

Now, observe that the second integral above goes to zero when $h \to 0$, as a consequence of equation \eqref{difW}. The first integral also goes to zero when $h \to 0$ 
 because, using the fact that $$\frac{e^{ah}-1}{h}-a=\sum_{k\geq 2} \frac{a^k h^{k-1}}{k!}\,,$$
 with $a=\frac{W(y|x+he_j)-W(y|x)}{h} $,
 we have
 \begin{eqnarray*}
 & & \bigg| \frac{e^{W(y|x+he_j)-W(y|x)}-1}{h} - \frac{W(y|x+he_j)-W(y|x)}{h} \bigg|=
\bigg| \sum_{k\geq 2}  \frac{\left(W(y|x+he_j)-W(y|x)\right)^k}{k!h}
\bigg| \\
 & & = \bigg| \sum_{k\geq 2}  \frac{\left(W(y|x+he_j)-W(y|x)\right)^{k-1}}{k!}
\frac{\left(W(y|x+he_j)-W(y|x)\right)}{h}
\bigg| \to 0 \;\;\mbox{when } \,h \to 0\;.
 \end{eqnarray*}

In the last expression we used the uniform continuity of $W$ and also \eqref{difW}.


\cqd

Remark: In the case where $A$ depends only on the two first  coordinates, we have that (see \eqref{op2coord})
$$\psi_A(x_1)=\frac{1}{\lambda_A}\int_{S^1}e^{A(y_1,x_1)}\psi_A(y_1)\,d\nu(y_1) .$$
Hence, $\psi_A$ satisfies the equation
$$\frac{\partial }{\partial x_1}\psi_A(x_1)= \frac{1}{\lambda_A}\int_{S^1}e^{A(y_1,x_1)}D_2A(y_1,x_1)\,\psi_A(y_1)\,d\nu(y_1).$$}

\end{document}